\documentclass[final,5p,times,twocolumn]{elsarticle}

\usepackage{amsmath,amssymb,bbm,comment,enumitem,hyperref,mathrsfs,algorithm,algorithmic,soul,comment,dsfont,xcolor}

\newtheorem{definition}{Definition}
\newtheorem{assumption}{Assumption}
\newtheorem{proposition}{Proposition}
\newtheorem{theorem}{Theorem}
\newtheorem{lemma}{Lemma}
\newtheorem{claim}{Claim}

\newtheorem{remark}{Remark}
\newtheorem{example}{Example}

\newcommand{\la}{\left\langle}
\newcommand{\ra}{\right\rangle}
\newcommand{\opt}{*}
\newcommand{\stat}{{\scriptscriptstyle \infty}}
\newcommand{\hc}{\mbox{\tiny hc}}

\newcommand{\blue}[1]{{\textcolor{black}{#1}}{}}

\journal{Operations Research Letters (June, 2021). Accepted (May, 2022)}

\begin{document}

\begin{frontmatter}


\title{Asymptotically Optimal Idling in the $GI$/$GI$/$N$+$GI$ Queue}
\author[label1]{Yueyang Zhong}
\author[label1]{Amy R. Ward}
\author[label2]{Amber L. Puha\fnref{cor1}}
\fntext[cor1]{Research supported by NSF Grant DMS-1712974}
\address[label1]{The University of Chicago Booth School of Business}
\address[label2]{California State University San Marcos}


\begin{abstract}
We formulate a control problem for a $GI$/$GI$/$N$+$GI$ queue, whose objective is to trade off the long-run average operational costs with server utilization costs. To solve the control problem, we consider an asymptotic regime in which the arrival rate and the number of servers grow large. The solution to an associated fluid control problem motivates that non-idling service disciplines are not in general optimal, unless some arrivals are turned away.  We propose an admission control policy designed to ensure that servers have sufficient idle time, which we show is asymptotically optimal. 
\end{abstract}


\end{frontmatter}


\section{Introduction} \label{section:intro}

One common assumption when studying the $GI$/$GI$/$N$+$GI$ queue is that the service discipline is non-idling; that is, that servers do not idle when customers are present in the queue (\cite{whitt2006fluid,kang2010fluid, kang2012asymptotic, zhang2013fluid, kang2019equivalence}).  However, in the restricted $M$/$M$/$N$+$M$ setting, the paper~\cite{ZW19} (see Theorem 1, Proposition 1, and Example 1 therein) shows that in the presence of server utilization costs, a non-idling service discipline may not be asymptotically optimal.  
\blue{Our purpose in this paper is to show that a similar phenomenon occurs in the $GI$/$GI$/$N$+$GI$ setting; that is, a non-idling service discipline might be suboptimal in the non-Markovian setting, when the system operates in a first-come, first-served (FCFS) manner}.

The $GI$/$GI$/$N$+$GI$ queue is more difficult to analyze than the $M$/$M$/$N$+$M$ queue because the state descriptor is more complex.  In particular, tracking the one-dimensional number-in-system process is sufficient when studying the $M$/$M$/$N$+$M$ queue, but more is needed when studying the $GI$/$GI$/$N$+$GI$ queue.  This is because a Markovian state descriptor must also include knowledge regarding the time that has elapsed since the last arrival, the amount of time each job in service has been in service, and the amount of time each job in the queue has waited, resulting in a measure-valued state descriptor.

The control question is to determine when an available server should take the next customer into service, and when such a server should idle for some period of time.  Too much idleness may lead to customer abandonment and excessive waiting, whereas too little rest may lead to server fatigue.  To quantify these two competing interests, we consider an objective function that trades off the abandonment costs (and also, as an extension, holding costs) with server utilization costs.  Exact analysis of the $GI$/$GI$/$N$+$GI$ queue is intractable, and, therefore, we study the queue in an overloaded asymptotic regime in which the arrival rate and the number of servers become large.  In that regime, we formulate a fluid control problem, and find that the solution to the fluid control problem sometimes motivates idling servers when customers are waiting (when operational costs are small compared to utilization costs).  The policy we propose, and show is asymptotically optimal (see our main results in Theorems~\ref{theorem:convergence-cost} and~\ref{theorem:lower-bound}, and their extension to incorporate holding costs in the online appendix), is one that ``thins'' the arrival process just enough to ensure the server utilization matches the solution to the fluid control problem.

Incorporating server utilization in the objective function is one way to ensure that the service discipline does not overwork servers. This can lead to increased employee retention, which can have performance benefits (discussed in \cite{wardimpact}).  Not overworking servers means ensuring sufficient idleness for all servers, an idea that arose earlier in papers that studied how to be fair to heterogeneous servers that can be grouped into statistically identical pools (see, e.g., \cite{ASS2011}, \cite{ward2013blind}), and how to exploit heterogeneous customers preferences so as to maximize revenue (see, e.g., \cite{afeche2016optimal},~\cite{maglaras2018optimal}).


{\bf Notation.}
We denote the set of integers endowed with the discrete topology by $\mathbb{Z}$, the set of non-negative integers by $\mathbb{Z}_{+}$, the set of positive integers by $\mathbb{N}$, the set of real numbers endowed with the Euclidean topology by $\mathbb{R}$, and the set of non-negative real numbers by $\mathbb{R}_{+}$.
For $F$, a cumulative distribution function (abbreviated c.d.f.\ henceforth) on ${\mathbb R}_+$ with density $f$, we write $\bar F=1-F$ and recall that the right edge of the support is given by
$x_r=\sup\{x\in{\mathbb R}_+ : \bar F(x)>0\}$ and the hazard function is $x\mapsto f(x)/\bar F(x)$ for $x\in[0,x_r)$.
For a measurable space $(S, \mathcal{F})$ and a measurable set $A \in \mathcal{F}$, $1_A$ is the indicator function of the set $A$, which is one when its argument is a member of the set $A$ and is zero otherwise. In addition, when $A$ is $S$, we use the shorthand notation $1$ to mean $1_S$.
For $H\in (0,\infty]$, let $\mathbf{M}[0,H)$ denote the set of finite, non-negative Borel measures on $[0,H)$ endowed with the topology of weak convergence. For a given $\eta\in\mathbf{M}[0,H)$ and a Borel measurable function $f: [0,H) \rightarrow \mathbb{R}_{+}$ that is integrable with respect to $\eta$, we write $\la f,\eta \ra = \int_{[0,H)} f(x)\eta(dx)$. The set $\mathbf{M}[0,H)$ endowed with the topology of weak convergence is a Polish space~(\cite{prokhorov1956convergence}).
We let $\mathbf{0}\in\mathbf{M}[0,H)$ be the measure such that $\la f, \mathbf{0} \ra=0$ for all Borel measurable functions $f: [0,H) \rightarrow \mathbb{R}_{+}$. Given $x\in[0,H)$, $\delta_x$ denotes the Dirac measure in $\mathbf{M}[0,H)$ such that for all Borel measurable functions $f: [0,H)\rightarrow \mathbb{R}_{+}$, $\la f, \delta_x\ra=f(x)$. Then let $\mathbf{M}_D[0,H)$ denote the subset of $\mathbf{M}[0,H)$ consisting of the measures $\eta\in\mathbf{M}[0,H)$ such that either $\eta=\mathbf{0}$ or $\eta$ can be represented as a sum of finitely many Dirac measures, that is, $\eta=\sum_{i=1}^{n}a_i \delta_{x_i}$, for some finite $n\in\mathbb{N}$, $(a_1,\dots,a_n)\in(0,\infty)^n$ and $(x_1,\dots,x_n)\in[0,H)^n$. 
Given a Polish space $\mathbb{S}$, we use $\mathbf{D}(\mathbb{S})$ to denote the set of $\mathbb{S}$ valued functions of $\mathbb{R}_{+}$ that are right continuous with finite lefts, endowed with the usual Skorokhod $J_1$-topology. Finally, we use $\Rightarrow$ to denote weak convergence and $\overset{d}{=}$ to denote equivalence in distribution.

\vspace{-0.1in}
\section{The Model and Admissible Policy Class} \label{section:model}
\vspace{-0.1in}
In this paper, we study a single-class many server queue with generally distributed inter-arrival, service, and patience times (i.e., a $GI$/$GI$/$N$+$GI$ queue) operating under a head-of-the-line (HL) control policy, that may or may not be non-idling. This is as specified in~\cite{puha2020fluid} specialized to a single customer class.  In particular, we consider the model specified in~\cite{kang2010fluid}, but with the non-idling condition~\cite[(2.30)]{kang2010fluid} removed.  Absent the non-idling condition, the system dynamics are not uniquely specified.  Hence, one must specify a control policy to determine when each customer in system will commence service.  Such control policies should satisfy natural conditions such as not using information about the future to make scheduling decisions.  In what follows, we describe the model and admissible policy class \blue{in brief}.  
We refer the interested reader to~\cite{puha2020fluid} for details.

{\bf The Model.}
Customers arrive according to a delayed renewal process $E$ with rate $\lambda\in\mathbb{R}_{+}$, each with a service time sampled from c.d.f.\ $G^s$ having finite mean $1/\mu\in(0,\infty)$, and a patience time \blue{(also known as reneging time)} sampled from a c.d.f.\ $G^r$ having finite mean $1/\theta\in(0,\infty)$. We denote the c.d.f. for the inter-arrival distribution associated with the renewal arrival as $G$. We assume $G$, $G^s$ and $G^r$ are absolutely continuous with density functions $g$, $g^s$ and $g^r$ respectively that have right edges of support $H$, $H^s$ and $H^r$ respectively and hazard function $h$, $h^s$ and $h^r$ respectively.  We assume that there exists $0\leq L^s< H^s$ such that $h^s$ is either bounded or lower-semicontinuous on $(L^s,H^s)$
and $h^r$ is bounded and continuous. Boundedness of $h^r$ implies that $H^r=\infty$.
Finally, we assume $G^r$ is strictly increasing with inverse function $(G^r)^{-1}$. The queue indexed by $N\in\mathbb{N}$ has $N$ identical servers and is defined on a fixed probability space $(\Omega,\mathcal{F}, \mathbb{P})$. For the remainder of this paper, we superscript all quantities that depend on $N$ by $N$, e.g., $G^N$, $g^N$, $H^N$, $\lambda^N$ and $E^N$ depend on $N$, but $G^s$ and $G^r$ do not vary with $N$.

\blue{Following the notation in Section 2.2 in~\cite{puha2020fluid}}, the state descriptor for the $N$-server queue is denoted by $y^{N} = (\alpha^{N}, x^{N}, \nu^{N}, \eta^{N})\in\mathbb{Y}_D$, where $\mathbb{Y}_D = \mathbb{R}_{+} \times \mathbb{Z}_{+} \times \mathbf{M}_D[0,H^s) \times \mathbf{M}_D[0,H^r)$. In particular, $\alpha^{N}\in[0,H^N)$ is the time that has elapsed since the last customer arrived to the system, $x^{N}\in\mathbb{Z}_{+}$ is the number of customers in system, $\nu^{N}\in\mathbf{M}_{D}[0,H^s)$ is a measure that has a unit mass at the age-in-service (amount of service received) of each customer currently in service, and $\eta^{N}\in\mathbf{M}_{D}[0,H^r)$ is a measure that has a unit mass at the potential waiting time of each customer ``potentially'' in system, (that is, each unit mass tracks the time passed since a customer's arrival, until that customer's patience time expires, at which point the unit atom is removed and tracking stops.) When $Y^N(0)$ denotes the initial state, the coordinate $\alpha^N(0)$ determines the distribution of the initial delay for $E^N$ as the conditional distribution of $G^N$ given $\alpha^N(0)$. That is, the initial delay distribution has density $g_0^N(x)=\frac{g^N(\alpha^N(0)+x)}{1-G^N(\alpha^N(0))}$ for $x \in [0,H^N-\alpha^N(0))$.

A state process for the $N$-server queue is a $\mathbb{Y}_D$ valued, right continuous process $Y^N$ with finite left limits that satisfies a set of dynamic equations for the $N$-server queue consistent with HL service. These are given as equations~(5)-(26) in~\cite{puha2020fluid}, which we omit here due to space constraints. With these, customers can only enter service at or after their arrival time and prior to their patience time expiring. An available server may idle or may take the customer in queue with the largest waiting time, the HL customer, into service. Once a server commences serving a customer, it works at rate one on the work associated with that customer until completely fulfilling that customer's service requirement, at which point the customer departs.

{\bf \blue{The Admissible Policy Class.}}
The admissible policy class consists of all policies that only allow customers to enter service at moments of a customer departure or arrival, do not use information about the future, 
and are such that the state process $Y^N$ is a Feller Markov process with respect to a natural filtration, and whose initial condition is policy compatible. The following leverages~{\cite{puha2020fluid}} to make this more precise. 

As mentioned above, equations~(5)-(26) in~\cite{puha2020fluid}
do not uniquely specify the system dynamics.
These are uniquely determined by the specification of an HL control policy $\pi^N=({\mathbb S}^N,\{{\mathbb P}_y^N\}_{y\in{\mathbb S}^N})$. Here, as in Definition~1 in~\cite{puha2020fluid}, ${\mathbb S}^N$ is the Polish subspace of ${\mathbb Y}_D$ that corresponds to the set of states that are achievable under the control policy.  Also, for each initial state $y\in{\mathbb S}^N$, ${\mathbb P}_y^N$ is 
a probability measure that uniquely determines the system dynamics when the system starts in state $y$.  More formally, $\{{\mathbb P}_y^N\}_{y\in{\mathbb S}^N}$ is a collection of probability measures indexed by ${\mathbb S}^N$ such that the mapping $y\mapsto {\mathbb P}_y^N(B)$ from ${\mathbb S}^N$ to $[0,1]$ is Borel measurable for each measurable $B\subset{\mathbf D}({\mathbb S}^N)$ and, for each $y\in{\mathbb S}^N$, ${\mathbb P}_y^N$ almost surely,
\vspace{-0.05in}
\begin{align}\label{eq:admissible-kernal}
	Y^N(0)=y,\ Y^N\in {\mathbf D}({\mathbb S}^N)\hbox{ and satisfies }(5)-(26)\hbox{ in~\cite{puha2020fluid}}.
\end{align}

\vspace{-0.05in}	
Given an HL control policy $\pi^N$, a state process $Y^N$ satisfying (\ref{eq:admissible-kernal}) specifies an entry-into-service process $K^N$.  Indeed, since a job has age-in-service equal to zero at the time of entering service, $\la 1_{\{0\}},\nu^N(t)\ra$ is the number of jobs to enter service at time $t$, for each $t>0$.  Then $K^N$ is a counting process such that $K^N(0)=0$ and $K^N(t)-K^N(t-)=\la 1_{\{0\}},\nu^N(t)\ra$ for each $t>0$.  In particular, $K^N(t)$ is the number of customers that enter service by time $t$ for each $t\ge 0$. Then, for each $t\ge 0$, $D^N(t)=\langle 1,\nu^N(0)\rangle+K^N(t)-\langle 1,\nu^N(t)\rangle$ denotes the number of customers to depart the system due to service completion by time $t$. We restrict attention to HL policies that only allow customers to enter service at moments of a customer departure or arrival.  We require that for each $y\in{\mathbb S}^N$, ${\mathbb P}_y^N$ almost surely, for all $t \geq 0$,
\vspace{-0.08in}
\begin{equation}\label{eq:KControl}
K^N(t)-K^N(t-)\le E^N(t)-E^N(t-)+D^N(t)-D^N(t-).\footnote{This condition is sufficient for a tightness result to hold as shown in \cite{puha2020fluid}.}
\end{equation}
 
\vspace{-0.08in}
We allow for random initial states that are compatible with a given HL control policy $\pi^N=({\mathbb S}^N,\{{\mathbb P}_y^N\}_{y\in{\mathbb S}^N})$.
As in Definition~2 in~\cite{puha2020fluid}, an initial distribution for $\pi^N$ is a Borel probability measure $\varsigma^N$ on ${\mathbb S}^N$ that determines the distribution of the initial state $Y^N(0)$.
In particular, for each measurable $B\subset{\mathbf D}({\mathbb S}^N)$, define $\mathbb{P}^{N}_{\varsigma}(B)=\int_{{\mathbb S}^N}{\mathbb P}_y^N(B)\varsigma^N(dy)$.
Then $\mathbb{P}^{N}_{\varsigma}$ denotes the distribution of the state process $Y^N$ under $\pi^N$ for initial distribution $\varsigma^N$.
We say that an initial distribution $\varsigma^N$ for $\pi^N$ is compatible if $\mathbb{E}^N_{\varsigma}\left[\la 1,\eta^N(0)\ra\right]<\infty$, where $\mathbb{E}^{N}_{\varsigma}$ denotes the expectation operator for $\mathbb{P}^{N}_{\varsigma}$. Given an HL control policy $\pi^N$ and a compatible initial distribution $\varsigma^N$, we refer to the process $Y^N$ with law $\mathbb{P}^N_{\varsigma}$ as the state process for $(\pi^N, \varsigma^N)$.

In order to restrict attention to HL control policies that do not use information about the future, we require $K^N$ to be non-anticipating. This amounts to requiring $K^N$ to be adapted to a suitable filtration as in Definition~3 in~\cite{puha2020fluid}.  
Because we consider long-run average cost, we make a further restriction in the definition of admissible HL control policies, which is used in Section~\ref{section:ao} to establish the existence of a stationary distribution.

\vspace{-0.1in}
\begin{definition}[Admissible Policies] \label{definition:admissible-policy}
An admissible HL control policy for $E^N$ is an HL control policy $\pi^N$ such that for any compatible initial distribution $\varsigma^N$, the pair $(\pi^N,\varsigma^N)$ (i) satisfies Definition 3 in~\cite{puha2020fluid} and \eqref{eq:KControl} and
(ii) is such that the state process $Y^N$ for $(\pi^N,\varsigma^N)$ is a Feller Markov process with respect to the filtration used in Definition 3 in~\cite{puha2020fluid}.	
\end{definition}

\vspace{-0.2in}
\begin{remark}
\blue{Our admissible policies focus on HL (equivalently, FCFS) control policies due to their common use in practice. However, non-HL control policies can be optimal in some settings; see~\cite{bassamboo2016scheduling}.}
\end{remark}

\vspace{-0.1in}
Let $\Pi^N$ denote the set of admissible HL control policies for $E^N$ in Definition~\ref{definition:admissible-policy}. For $\pi^N\in\Pi^N$, we will sometimes write $Y^N(\pi^N,\cdot)$, $X^N(\pi^N,\cdot)$, $\nu^N(\pi^N,\cdot)$, $\eta^N(\pi^N,\cdot)$, $K^N(\pi^N,\cdot)$ or $D^N(\pi^N,\cdot)$ to make the dependence on $\pi^N$ explicit.

\begin{proposition} \label{prop:stationary-exist}
For any $\pi^N\in\Pi^N$, there exists a compatible initial distribution $\xi^N$ such that the state process $Y_\stat^N$ for $(\pi^N, \xi^N)$ is a stationary process.
\end{proposition}
Proposition~\ref{prop:stationary-exist} follows as a special case of Lemma~\ref{lemma:stationary-exist} stated in Section~\ref{section:prelim-results}.

Given $\pi^N\in\Pi^N$ and a compatible initial distribution $\xi^N$ such that the state process $Y_\stat^N$ for $(\pi^N, \xi^N)$ is a stationary process, we refer to $\xi^N$ as a compatible stationary distribution for $\pi^N$
and we let ${\mathcal S}(\pi^N)$ denote the set of all compatible stationary distributions for $\pi^N$.


\section{The Control Problem} \label{section:scheduling}
Each customer abandonment incurs a cost $a\in(0,\infty)$ and the strictly increasing, continuous and convex function $g_U: [0,1] \rightarrow[0,\infty)$ captures the cost of server utilization. The trade-off is between working the servers as much as possible, which incurs high utilization cost but low abandonment cost, and giving the servers more rest, which incurs lower utilization cost but higher abandonment cost. 
In particular, given $\pi^{N}\in\Pi^N$ and a compatible initial distribution $\varsigma^N$, we define the long-run average cost of $(\pi^N,\varsigma^N)$ as
\begin{align*}
\mathcal{C}_{\varsigma}^{N}(\pi^{N}) := \limsup_{T\rightarrow\infty}
 \frac{1}{T} \mathbb{E}^{N}_{\varsigma} \left[a \frac{R^{N}(\pi^{N},T)}{N} + \int_{0}^{T} g_U\left(\frac{B^{N}(\pi^{N},t)}{N}\right)dt\right],
\end{align*}
where, for each $t > 0$, $R^{N}(\pi^N,t)$ is the cumulative number of abandonments by time $t$ under $\pi^N$, and $B^{N}(\pi^N,t) \leq N$ is the number of busy servers at time $t$ under $\pi^N$.

\begin{proposition} \label{prop:stationary-dist}
	For any $\pi^N\in\Pi^N$ and compatible initial distribution $\varsigma^N$, there exists $\xi^N\in{\mathcal S}(\pi^N)$ such that $\mathcal{C}_{\varsigma}^N(\pi^N) = \mathcal{C}_{\xi}^N(\pi^N)$.
\end{proposition}
Proposition~\ref{prop:stationary-dist} follows as a special case of Lemma~\ref{lemma:stationary-dist} stated in Section~\ref{section:prelim-results}.

Given $\pi^N\in\Pi^N$, let 
$
\mathcal{C}^{N}(\pi^{N}) :=\sup_{\xi^N\in{\mathcal S}(\pi^N)}\mathcal{C}_{\xi}^{N}(\pi^{N})
$
denote the worst case cost.
By Proposition \ref{prop:stationary-dist}, $\mathcal{C}^{N}(\pi^{N})$ is the supremum of $\mathcal{C}_{\varsigma}^{N}(\pi^{N})$ over all compatible initial distributions $\varsigma^N$.
Our objective is to find an admissible control policy $\pi^{N}_{\mbox{\tiny opt}}$ such that 
\begin{align} \label{eq:scheduling-original}
\mathcal{C}^{N}(\pi^{N}_{\mbox{\tiny opt}}) := \inf_{\pi^{N}\in\Pi^{N}}\mathcal{C}^{N}(\pi^{N}).
\end{align}
The objective is such that a non-idling control policy is not in general optimal.
Based on the discrete-event queuing model, it is not possible to solve for $\pi^{N}_{\mbox{\tiny opt}}$ exactly. Thus, we leverage an analytically tractable approximating fluid control problem to postulate an HL control policy that one might expect to perform well for the objective~(\ref{eq:scheduling-original}). Then, we show that this policy is asymptotically optimal (see Theorems~\ref{theorem:convergence-cost} and~\ref{theorem:lower-bound} in Section~\ref{section:ao}).

\section{The Fluid Control Problem} \label{section:fluid}
The fluid control problem is based on the fluid model and the fluid model solutions defined in~\cite{puha2020fluid}.  
Fluid model solutions arise as functional law of large numbers limits of sequences of state descriptors for the stochastic system under fluid scaling. For each $N\in\mathbb{N}$, we define the fluid scaling for the $N$-server system as follows.  Recall the constant $\lambda^N$ and the processes $E^N$, $\alpha^N$, $X^N$, $\nu^N$, $\eta^N$, $K^N$ and $D^N$ defined in Section~\ref{section:model}, and the processes $R^N$ and $B^N$ defined in Section~\ref{section:scheduling}; also define the process $Q^N = X^N-B^N$ as the queue length, and the process $I^N=N-B^N$ as the number of idle servers. Then, let $\bar{\alpha}^N = \alpha^N$; also for $\triangle^N = \lambda^N, E^N, X^N, \nu^N, \eta^N, K^N, D^N, R^N, B^N$, $Q^N, I^N$, let $\bar{\triangle}^N = \triangle^N / N$. Then, the fluid-scaled state process for the $N$-server system is $\bar{Y}^N = (\bar{\alpha}^N, \bar{X}^N, \bar{\nu}^N, \bar{\eta}^N)$. Under suitable asymptotic conditions, limit points exist and are fluid model solutions almost surely (see Lemma \ref{lemma:weak-convergence} in Section \ref{section:prelim-results}.)

In particular, fluid model solutions are functions of time that take values in the set $\mathbb{X}=\mathbb{R}_{+} \times \mathbf{M}[0,H^s) \times \mathbf{M}[0,H^r)$ endowed with the product topology.  Then a state $(x,\nu,\eta)\in\mathbb{X}$ for the fluid model is a fluid analog of the state descriptor for the stochastic system with $x$, $\la 1_{[0,z]},\nu \ra$ and $\la 1_{[0,z]}, \eta \ra$ corresponding to the total mass in system, the total mass in service with age-in-service less than or equal to $z$ for each $z\in\mathbb{R}_{+}$, and the total mass potentially in system of age less than or equal to $z$ for each $z\in\mathbb{R}_{+}$, respectively.
They satisfy a set of conditions determined by a positive constant $\gamma$, which is the rate at which ``fluid" or mass arrives to the system.  These conditions are referred to as the fluid model for $\gamma$. We summarize the fluid model for $\gamma$ and the definition of a fluid model solution for $\gamma$ in~\ref{section:appendix}.

The invariant states for the fluid model for $\gamma$ are fixed points of the fluid model for $\gamma$. 
From Proposition 1 in~\cite{puha2020fluid}, an invariant state for $\gamma$ is determined by the long-run average fraction of the collective server effort provided to the customers, denoted by $b$. It is clear that $b$ must satisfy $b\in[0,\min\{1,\gamma/\mu\}]$, where we recall that $\mu$ is the reciprocal of the mean of $G^s$.
Then, when the initial state for a fluid model solution for $\gamma$ is an invariant state for $\gamma$, it turns out that the departure rate of the fluid from the system is $b\mu$ and so, by conservation of mass, $\gamma-b\mu$ must be the rate at which fluid abandons. This implies that the abandonment rate is insensitive to the patience time distribution, which has a similar flavor to the insensitivity result for a single server queue in the large deviations regime in~\cite{atar2019large}.

\begin{assumption}\label{assumption:lambda-asymptotics}
	Let $\lambda\in(0,\infty)$. Suppose that $\lim_{N\rightarrow\infty}\bar{\lambda}^N=\lambda$.
\end{assumption}
Henceforth, $\lambda$ satisfying the conditions in Assumption~\ref{assumption:lambda-asymptotics} is fixed. Our fluid control problem is based on the invariant states for $\lambda$.
We expect to obtain the following fluid control problem for $\lambda$ when letting $N\rightarrow\infty$ in problem~\eqref{eq:scheduling-original}.

\begin{definition}[The Fluid Control Problem]
The fluid control problem for $\lambda$ is given by
	\begin{align} \label{eq:scheduling-fluid}
		\min_{b\in[0,\min\{1,\lambda/\mu\}]} \ a(\lambda-b\mu) + g_U(b).
	\end{align}
\end{definition}
We denote the solution to~\eqref{eq:scheduling-fluid} by $b_{\opt}$ (which exists and is unique because~\eqref{eq:scheduling-fluid} optimizes a convex function over a compact set).

\begin{example}\label{example1}
Suppose $a=1$ and $g_U(b) = b^2$. Then, the solution to (\ref{eq:scheduling-fluid}) is $b_{\opt} = \min\{1, \mu/2, \lambda/\mu\}$.
\end{example}

The solution to~\eqref{eq:scheduling-fluid} motivates a control policy that we expect to have good performance with respect to the original objective~\eqref{eq:scheduling-original} when the arrival rate $\lambda^N$ and the number of servers $N$ are large.  When $b_{\opt}=\min\{1,\lambda/\mu\}$, we expect a non-idling control policy to be optimal for~\eqref{eq:scheduling-original}. Otherwise, when $b_{\opt}<\min\{1,\lambda/\mu\}$, the solution to the fluid control motivates defining a policy that uses customer abandonments to trim congestion, in order to reduce server workload, and provide (additional) server idle time. In this case, for each $N\in\mathbb{N}$, consider the HL control policy $\tilde{\pi}^N$ such that each server idles after each service completion for the difference between the desired expected time between service completions, $(b_{\opt}\mu)^{-1}$, and the expected time between service completions when the server is always busy, $\mu^{-1}$; that is, for $(b_{\opt}\mu)^{-1} – \mu^{-1} = (1-b_{\opt})(b_{\opt}\mu)^{-1}$ time units. Such a policy seems quite reasonable, and should be asymptotically optimal. However, establishing that for any sequence of compatible initial distributions $\{\varsigma^N\}_{N\in{\mathbb N}}$,
\begin{align}\label{eq:intention}
\begin{aligned}
	&\lim_{N \rightarrow \infty} \lim_{t \rightarrow \infty} \frac{1}{t}{\mathbb E}_{\varsigma}^N\left[\bar{R}^N(\tilde{\pi}^N,t)\right] = \lambda-b_{\opt}\mu \ 
	\mbox{ and }  \\
	&\lim_{N \rightarrow \infty} \lim_{t \rightarrow \infty} {\mathbb E}_{\varsigma}^N\left[g_U\left(\bar{B}^N(\tilde{\pi}^N,t)\right)\right] = g_U(b_{\opt})
\end{aligned}
\end{align}
is difficult. This difficulty is related to a lack of results providing sufficient conditions for fluid model solutions to converge to invariant states in the time infinity limit (see Section 7.1 in~\cite{kang2012asymptotic}).   Instead, we propose to expand the admissible policy class to include thinned arrival processes and then rely on results in the literature for non-idling many server queues to show that (\ref{eq:intention}) holds. \blue{If we can show a policy is asymptotically optimal  for an enlarged policy class, then we know that no policy in the original smaller policy class can perform better.} 

\vspace{-0.1in}
\section{The Proposed Policy $\pi^N_{\opt}$} \label{section:policy}
\vspace{-0.1in}
The solution $0 \le b_{\opt} \le \min\{1,\lambda/\mu\}$ to~\eqref{eq:scheduling-fluid} represents the optimal long-run average fraction of busy servers, which suggests that a control policy that thins the arrival process to rate $b_{\opt}\mu$ and forces the servers to work in a non-idling fashion, but builds in idleness due to admission control, should perform well for the original objective~\eqref{eq:scheduling-original}. This motivates us to enlarge the admissible policy class in Definition~\ref{definition:admissible-policy} to allow for admission control.  Specifically, at the time of each arrival, let $p \in (0,1]$ be the probability the arrival is admitted for service and $1-p$ the probability the arrival is rejected, which incurs a cost $a$. Given $p\in(0,1]$, we denote the admitted arrival process by $E_p^N$,  and we refer to the $N$-server queue with arrival process $E_p^N$ as the $p$-admitted queue. It is clear that the thinned arrival process $E_{p}^N$ is a suitably delayed renewal process with arrival rate $p\lambda^N$, because the admitted arrivals remain i.i.d..

\begin{definition}[Enlarged Admissible Policies] \label{definition:admissible-policy-2}
For any $p\in(0,1]$, an admissible HL control policy for $E_p^N$ satisfies Definition~\ref{definition:admissible-policy} with $E^N$ replaced by $E_p^N$.
\end{definition}
For $p\in(0,1]$, let $\Pi^N_p$ denote the set of admissible HL control polices for $E_p^N$. Note that $\Pi^N_1=\Pi^N$.  
For $p\in(0,1]$, $\pi^N_p \in \Pi^N_p$ and $\triangle^N =Y^N, X^N, \nu^N, \eta^N, K^N$, $D^N, R^N, B^N, Q^N$ or $I^N$, $\triangle^N(\pi^N_p,\cdot)$ refers to the process for the $p$-admitted queue under $\pi^N_p$.

Given $p\in(0,1]$, $\pi^N_p\in\Pi^N_p$ and a compatible initial distribution $\varsigma^N$, the long-run average cost of $(\pi^N_p,\varsigma^N)$ is
\begin{align}\label{eq:scheduling-2}
	\mathcal{C}_{\varsigma}^{N}(\pi^{N}_p) := \limsup_{T\rightarrow\infty}\frac{1}{T} \mathbb{E}^{N}_{\varsigma} &\left[a \left(\bar{E}^N(T)-\bar{E}^N_p(T)+\bar{R}^{N}(\pi^{N}_p,T)\right) \right. \nonumber\\
	& + \left.\int_{0}^{T} g_U\left(\bar{B}^{N}(\pi^{N}_p,t)\right)dt\right].
\end{align}
When the initial state for the fluid model for $p\lambda$ is an invariant state for $p\lambda$ associated with $b\in[0,p\lambda/\mu]$, $p\lambda-b\mu$ is the rate at which fluid abandons and $(1-p)\lambda$ is the rate at which fluid is rejected. Since $p\in(0,1]$ is a parameter that can be optimized over, the resulting fluid control problem is given by
\begin{align} \label{eq:scheduling-fluid-2}
	&\min_{p\in(0,1],b\in[0,\min\{1,p\lambda/\mu\}]}\  a\left(1-p\right)\lambda + a \left(p\lambda - b\mu\right) + g_U(b)\nonumber\\
	=&\min_{b\in[0,\min\{1,\lambda/\mu\}]}a\left(\lambda-b\mu\right) + g_U(b).
\end{align}
The solution to~(\ref{eq:scheduling-fluid-2}) does not depend on the admission control parameter $p\in(0,1]$ and is identical to the solution to~(\ref{eq:scheduling-fluid}). \blue{This observation crucially relies on the abandonment cost being linear with the per unit cost equal to the per unit cost of rejection}.

This gives us flexibility to propose a policy in $\Pi^N_p$ for various choices of $p\in(0,1]$. We first observe that an optimal admission control parameter must lie in $[b_{\opt}\mu/\lambda,1]$, because otherwise the admitted arrivals would not be sufficient for servers to work at busyness level $b_{\opt}$. Let
\begin{align}\label{eq:p_opt}
	p_{\opt}:=b_{\opt}\mu/\lambda.
\end{align}
We next observe that if the $p_{\opt}$-admitted queue satisfies the non-idling condition (that is, the servers never idle when customers are waiting), the long-run average fraction of busy servers achieves $b_{\opt}$. 
The non-idling condition, together with (5)-(26) in~\cite{puha2020fluid} uniquely specifies $\mathbb{P}^N_y$ for each $y\in\mathbb{S}^N=\left\{y^N\in\mathbb{Y}^D: N-\la 1,\nu^N \ra = (N-x^N)^{+} \mbox{ and } x^N \leq \la 1, \eta^N\ra\right\}$ and satisfies (\ref{eq:KControl}). Moreover, for any compatible initial distribution, the state process that satisfies the non-idling condition is a Feller, strong Markov process (see Proposition~4.2 in~\cite{kang2012asymptotic}). Thus, for any $p\in(0,1]$, the non-idling policy (the control policy that obeys the non-idling condition) is an admissible HL control policy for $E^N_p$, and thus is in $\Pi^N_p$.

\begin{definition}[The Proposed Policy] \label{definition:policy}
	For each $N\in\mathbb{N}$, let $\pi^N_{\opt}$ be the non-idling policy in $\Pi^N_{p_{\opt}}$, where $p_{\opt}$ is given by (\ref{eq:p_opt}).
\end{definition}


\section{Asymptotic Optimality of $\pi^N_{\opt}$} \label{section:ao}
In this section, we state our main results concerning asymptotic optimality of $\{\pi^N_{\opt}\}_{N\in\mathbb{N}}$ under fluid scaling.

\begin{theorem}[Convergence under the Proposed Policy]\label{theorem:convergence-cost} 
	Suppose that Assumption~\ref{assumption:lambda-asymptotics} holds \blue{and that $h^s$ is non-increasing when $b_\opt=1$}. Then the sequence $\{\pi^N_{\opt}\}_{N\in\mathbb{N}}$ satisfies
	\begin{align*}
		\lim_{N\rightarrow\infty} \mathcal{C}^N\left(\pi^N_{\opt}\right) = a(\lambda - b_{\opt} \mu) + g_U(b_{\opt}).
	\end{align*}
\end{theorem}

Let $\hat{\Pi}^N:=\cup_{p\in(0,1]} \Pi^N_p$ denote the enlarged policy class, and given $\hat{\pi}^N\in\hat{\Pi}^N$, let $\hat{p}^N\in(0,1]$ denote the associated admission control parameter.
\begin{theorem}[Asymptotic Lower Bound] \label{theorem:lower-bound}
Suppose that Assumption~\ref{assumption:lambda-asymptotics} holds, $\hat{\pi}^{N}\in\hat\Pi^N$ for each $N\in{\mathbb N}$ and the sequence 
$\{\hat{p}^N\}_{N\in\mathbb{N}}$ satisfies $\lim_{N\rightarrow\infty}\hat{p}^N = p$ for some $p\in(0,1]$. Then, 
\begin{align*}
	\liminf_{N\rightarrow\infty} \mathcal{C}^N(\hat{\pi}^{N}) \geq a(\lambda-b_{\opt} \mu) + g_U(b_{\opt}).
\end{align*}
\end{theorem}
\begin{remark}\label{remark:pN}
The condition that $\lim_{N\rightarrow\infty}\hat{p}^N = p$ for some $p\in(0,1]$ implies that $\{\hat{p}^N\bar{\lambda}^N\}_{N\in\mathbb{N}}$ satisfies $\lim_{N\to\infty}\hat{p}^N\bar{\lambda}^N=p\lambda$.
\end{remark}

Theorem~\ref{theorem:convergence-cost} establishes that the solution to the fluid control problem~\eqref{eq:scheduling-fluid} is achieved in the limiting system, when, for each $N$, the $N$-server system operates under $\pi^N_{\opt}$ in Definition~\ref{definition:policy}, and in case $b_\opt=1$, $h^s$ is non-increasing. Theorem~\ref{theorem:lower-bound} establishes that the fluid control problem~\eqref{eq:scheduling-fluid} is an asymptotic lower bound for the objective~\eqref{eq:scheduling-2}. As a consequence, we conclude that the proposed sequence of policies $\{\pi^N_{\opt}\}_{N\in\mathbb{N}}$ is asymptotically optimal.

The proof of Theorem~\ref{theorem:convergence-cost} given in Section~\ref{section:proofs-ao} is facilitated by the fact that, for each $N\in\mathbb{N}$, under $\pi_*^N$ the $p_\opt$-admitted $N$-server queue is non-idling, and thus, we can appeal to results in~\cite{kang2012asymptotic, atar2021large} to establish the weak convergence of the sequence of fluid-scaled stationary distributions. The additional condition that $h^s$ is non-increasing when $b_\opt=1$, is needed for this in order to apply part (3) of Theorem 3.2 in~\cite{atar2021large} in that case. This implies that the limit is the unique invariant state with zero queue mass.

The proof of Theorem~\ref{theorem:lower-bound} in Section~\ref{section:proofs-ao} requires first adapting one of the arguments in~\cite{kang2012asymptotic} (wherein the non-idling condition is assumed throughout) to show that a sequence of fluid-scaled stationary distributions is tight, and second arguing that the fluid control problem~\eqref{eq:scheduling-fluid-2} provides an asymptotic lower bound on the cost along any convergent subsequence.

In the next section, we establish some preliminary results for stationary distributions (for both the stochastic $N$-server queue model and the fluid model) that help to prove Theorems~\ref{theorem:convergence-cost} and~\ref{theorem:lower-bound}, which may also be of independent interest. The proofs of Theorems~\ref{theorem:convergence-cost} and~\ref{theorem:lower-bound} will be provided in Section~\ref{section:proofs-ao}.

\vspace{-0.08in}
\section{Preliminary Results}\label{section:prelim-results}
\vspace{-0.04in}
In order to prove our main results (Theorems~\ref{theorem:convergence-cost} and~\ref{theorem:lower-bound}), we begin by establishing two foundational results concerning stationary distributions for the $N$-server queue. Then, we provide a fluid limit theorem, which shows that the distributional limit points of stationary distributions are fluid model solutions almost surely under suitable asymptotic conditions. Finally, we show some properties of stationary fluid model solutions for $\gamma$. The proofs are delayed to the online appendix~\ref{section:proofs}.

{\bf Stationary Distributions of the $N$-Server Queue. }
The following lemmas confirm the existence of a stationary distribution under any admissible HL control policy for $E_p^N$ and $p\in(0,1]$, and derive an expression for the long-run average cost. 
We denote by $\triangle^N_{\stat}$ a stationary process associated with the process $\triangle^N$, for $\triangle^N = E^N, Y^N, X^N, \nu^N, \eta^N, K^N, D^N, R^N, B^N, Q^N, I^N$.

\begin{lemma} \label{lemma:stationary-exist}
Let $p\in(0,1]$. For any $\pi^N_{p}\in\Pi^N_p$, there exists a compatible initial distribution $\xi^N$ such that the state process $Y^N_\stat$ for $(\pi^N_p, \xi^N)$ is stationary. Moreover, $\mathbb{E}^N_{\xi}\left[ \la 1, \eta^N_\stat(t) \ra\right] = p\lambda^N\theta^{-1} <\infty$, for all $t\geq 0$.
\end{lemma}

\begin{remark}
Proposition~\ref{prop:stationary-exist} in Section~\ref{section:model} follows by setting $p=1$.
\end{remark}

Given $p\in(0,1]$, $\pi^N_p\in\Pi^N_p$ and a compatible initial distribution $\varsigma^N$, let 
\begin{align} \label{eq:def-chi-N}
	\chi^N(t):=\inf\{x\geq 0: \la 1_{[0,x]}, \eta^N(t)\ra \geq Q^N(t)\}
\end{align}
represent the waiting time of the HL customer at time $t$ for each $t\geq 0$. Then, for $t \geq 0$,
\begin{align}\label{eq:Q-chi}
	Q^N(t) = \la 1_{[0,\chi^N(t)]}, \eta^N(t) \ra.
\end{align}
The associated stationary process is denoted by $\chi^N_\stat$.

\begin{lemma}\label{lemma:stationary-dist}
Let $p\in(0,1]$. For any $\pi^N_{p}\in\Pi^N_p$ and compatible initial distribution $\varsigma^N$, there exists $\xi^N \in \mathcal{S}(\pi^N_p)$ such that
\begin{align}\label{eq:long-run-avg-R} 
	\resizebox{.86\hsize}{!}{$\limsup\limits_{T\to\infty} \mathbb{E}^N_{\varsigma}\left[\frac{\bar{R}^N(\pi^N_p,T)}{T}\right] = \mathbb{E}^N_{\xi}\left[\la 1_{[0,\chi^N_{\stat}(0)]} h^r, \bar{\eta}^N_{\stat}(0) \ra\right],$}
\end{align}
and
\begin{align}\label{eq:long-run-avg-B}
	\resizebox{.89\hsize}{!}{$\limsup\limits_{T\to\infty} \mathbb{E}^N_{\varsigma}\left[\frac{1}{T}\int_{0}^{T}g_U\left(\bar{B}^N(\pi^N_p,t)\right)dt\right] = \mathbb{E}^N_{\xi}\left[g_U\left(\bar{B}^N_{\stat}(0)\right)\right].$}
\end{align}
If $\varsigma^N \in \mathcal{S}(\pi^N_p)$, then $\xi^N=\varsigma^N$.
\end{lemma}
In light of~(\ref{eq:Q-chi}), one can interpret the right-hand side of~(\ref{eq:long-run-avg-R}) as an expected stationary reneging rate for the $N$-server queue. 
\begin{remark}\label{remark:stationary-dist}
For any $p\in(0,1]$, $\pi^N_p \in \Pi^N_p$ and compatible initial distribution $\varsigma^N$, there exists $\xi^N\in\mathcal{S}(\pi^N_p)$ such that 
\begin{align*}
	\mathcal{C}_{\varsigma}^N(\pi^N_p) = \mathbb{E}^N_{\xi} \left[a (1-p)\bar{\lambda}^N + a \la 1_{[0,\chi^N_{\stat}(0)]} h^r, \bar{\eta}^N_{\stat}(0) \ra + g_U\left( \bar{B}^N_\stat(0) \right)\right].
\end{align*}
Proposition~\ref{prop:stationary-dist} in Section~\ref{section:model} follows by setting $p=1$.
\end{remark}

{\bf A Fluid Limit Theorem. }
Here we provide asymptotic assumptions under which it is shown in~\cite{puha2020fluid} that fluid limit points are almost surely fluid model solutions. Such a result is crucial for the proof of Theorem~\ref{theorem:lower-bound}, which will appear in Section~\ref{section:proofs-ao}.

\vspace{-0.08in}
\begin{assumption}\label{assumption:asymptotics}
Suppose for each $N\in{\mathbb N}$, $p^N\in(0,1]$, $\pi^N_{p^N}\in\Pi^N_{p^N}$ for $E^N_{p^N}$ and $\varsigma^N$ is a compatible initial distribution for $\pi^N_{p^N}$. Assume that $\lim_{N\to\infty}p^N=p$ and $(\bar X^N(0), \bar \nu^N(0), \bar \eta^N(0)) \Rightarrow (X^0, \nu^0, \eta^0)$, as $N\rightarrow\infty$, for some random variable $(X^0, \nu^0, \eta^0)$ taking values in $\mathbb{X}$ such that $\sup_{N\in\mathbb{N}}\mathbb{E}^N_{\varsigma}\left[\la 1, \bar{\eta}^N(0) \ra\right]<\infty$.
\end{assumption}
\vspace{-0.1in}
\begin{remark}\label{remark:assumptions}
Under Assumptions~\ref{assumption:lambda-asymptotics} and~\ref{assumption:asymptotics} and the conditions on $E^N_{p^N}$, $K^N$, $G^s$, $g^s$, $h^s$, $G^r$, $g^r$, and $h^r$ specified in Sections~\ref{section:model} and~\ref{section:policy}, one can without loss of generality assume that the convergence of the initial condition in Assumption~\ref{assumption:asymptotics} is almost sure and then check that Assumptions~1, 2, 3(1), 3(3), 3(4), 4, 5(1) and~5(3) in~\cite{puha2020fluid} hold, i.e., Assumptions~3(2), 3(5) and 5(2) may not hold.
\end{remark}

\vspace{-0.1in}
\begin{lemma}\label{lemma:convergence-eta}
Suppose Assumptions~\ref{assumption:lambda-asymptotics} and~\ref{assumption:asymptotics} hold. Then, $\bar{\eta}^N \Rightarrow \eta$, as $N\to\infty$, where $\eta(0)\overset{d}{=}\eta^0$ and $\eta$ satisfies~(\ref{eq:fluid-47}) almost surely for $E(t)=p\lambda t$, $t\geq 0$. 
\end{lemma}

\vspace{-0.05in}
In fact, Assumptions~3(2) and 3(5) in~\cite{puha2020fluid} can be replaced by the condition $\sup_{N\in\mathbb{N}} \mathbb{E}^N_{\xi}\left[\la 1, \bar{\eta}^N(0) \ra\right] < \infty$ and Assumption~5(2) ($\eta^0$ has no atoms) is used to establish convergence of the scaled reneging processes to the expression in~(\ref{eq:fluid-40}. Thus, the result in Theorem~1 in~\cite{puha2020fluid} continues to hold. We obtain the following slightly restated version of Theorem~1 in~\cite{puha2020fluid}. 

\vspace{-0.05in}
\begin{lemma}[Theorem 1 in~\cite{puha2020fluid}] \label{lemma:weak-convergence}
Suppose that $\{(\pi^N,\varsigma^N)\}_{N\in\mathbb{N}}$ is such that Assumptions~\ref{assumption:lambda-asymptotics} and~\ref{assumption:asymptotics} hold, $\eta^0$ has no atoms, and $(X,\nu,\eta)$ is a distributional limit point of $\left\{ \left(\bar X^N, \bar \nu^N, \bar \eta^N  \right) \right\}_{N \in \mathbb{N}}$. Then $(X(0),\nu(0),\eta(0)) \overset{d}{=} (X^0, \nu^0, \eta^0)$ and $(X,\nu,\eta)$ is almost surely a fluid model solution for $p\lambda$.
\end{lemma}

{\bf Properties of Stationary Fluid Model Solutions. }
Fix $\gamma>0$. Here we consider the fluid model for $\gamma$ with random initial states such that the resulting fluid model solution is a stationary process. Lemmas~\ref{lemma:etastat} and~\ref{lemma:thm1} below, provide properties of such solutions. The proof of Theorem~\ref{theorem:lower-bound} relies on Lemmas~\ref{lemma:etastat} and~\ref{lemma:thm1}.

In what follows, we fix a fluid model solution $Z_\stat = (X_\stat, \nu_\stat, \eta_\stat)$ for $\gamma$ such that $Z_\stat$ is a stationary process. We denote the law of $Z_{\stat}(0)$ by $\xi$ and the expectation operator by $\mathbb{E}_{\xi}$. In addition, we define a Borel probability measure $\eta_e$ satisfying $d\eta_e(x)=\theta\bar G^r(x)dx$ for all $x\in\mathbb{R}_{+}$, where the subscript $e$ is mnemonic for excess life distribution. 

\begin{lemma}\label{lemma:etastat}
For all $t \geq 0$, $\eta_\stat(t)=\gamma\theta^{-1}\eta_e$. In particular, for all $t\geq 0$, $\eta_\stat(t)$ has no atoms, $x\mapsto\la 1_{[0,x]}, \eta_\stat(t) \ra$ is a continuous strictly increasing function on ${\mathbb R}_+$, and $\la 1,\eta_\stat(t)\ra=\gamma\theta^{-1}$.
\end{lemma}

\begin{lemma}\label{lemma:thm1}
There exists $b\in[0,\min\{1,\gamma/\mu\}]$ such that for all $t \geq 0$, $\mathbb{E}_{\xi} \left[B_\stat(t)\right] =b$ and $\mathbb{E}_{\xi}\left[\la  1_{[0,\chi_\stat(t)]}h^r, \eta_\stat(t) \ra\right] = \gamma - b\mu$.
\end{lemma}

\section{Proofs of Main Results (Theorems~\ref{theorem:convergence-cost} and~\ref{theorem:lower-bound})\label{section:proofs-ao}}
\textit{Proof of Theorem~\ref{theorem:convergence-cost}. }
For each $N \in \mathbb{N}$, let $\xi^N\in{\mathcal S}(\pi_*^N)$ which exists by Lemma~\ref{lemma:stationary-exist}, and recall that $Y^N_\stat(0)$ has distribution $\xi^N$. Consider the sequence $\{(\bar{X}_\stat^N(0), \bar{\nu}_\stat^N(0), \bar{\eta}_\stat^N(0))\}_{N\in\mathbb{N}}$. We wish to show that $\lim_{N\rightarrow\infty} \mathcal{C}(\pi^N_\opt)=a(\lambda-b_\opt \mu)+g_U(b_\opt)$. By Lemma~\ref{lemma:stationary-dist}, it suffices to show that, 
\begin{align} \label{eq:thm1-framework}
	&\lim_{N\to\infty}\mathbb{E}^N_\xi \left[a (1-p_\opt)\bar{\lambda}^N+ a \la 1_{[0,\chi^N_\stat(0)]}h^r, \bar{\eta}^N_\stat(0) \ra + g_U(\bar{B}^N_\stat(0))\right]\nonumber\\
	&\qquad= a(\lambda-b_\opt \mu)+g_U(b_\opt).
\end{align}
Note that $p_{\opt} \bar{\lambda}^N \rightarrow p_{\opt} \lambda$, as $N\rightarrow\infty$ (from Assumption~\ref{assumption:lambda-asymptotics}). This, together with the assumptions on $E^N$ (which $E^N_{p_{\opt}}$ inherits), $G^s$, $g^s$, $h^s$, $G^r$, $g^r$, and $h^r$ given in Section~\ref{section:model}, implies that Assumptions~3.1-3.5 in~\cite{kang2012asymptotic} hold for $\{(E_{p_*}^N,\pi_*^N,\xi^N)\}_{N\in{\mathbb N}}$. In addition, since it is assumed that $h^s$ is non-increasing when $b_{\opt}=1$, the result in Theorem~3.3 in~\cite{kang2012asymptotic} holds\footnote{There is a gap in the original proof of Theorem~3.3 in~\cite{kang2012asymptotic}, where a stationary distribution for the fluid model is assumed to coincide with the invariant state, which is unique since $G^r$ is strictly increasing. Under the conditions of Theorem~3.3 in~\cite{kang2012asymptotic}, Theorem~3.2(1) in~\cite{atar2021large} implies that this is true when  $b_{\opt}<1$. With the added condition that $h^s$ is non-increasing, Theorem~3.2(3) in~\cite{atar2021large} implies that this is true when $b_{\opt}=1$. Hence, the result in Theorem~3.3 in~\cite{kang2012asymptotic} holds in the present setting. See the discussion in~\cite{atar2021large} that follows the statement of Theorem~3.2 for a detailed explanation.}, which establishes
\begin{align} \label{eq:thm1-aux0}
	(\bar{X}_\stat^N(0), \bar{\nu}_\stat^N(0), \bar{\eta}_\stat^N(0)) \Rightarrow (b_\opt, b_\opt \nu_e, p_*\lambda\theta^{-1} \eta_e), 
\end{align}
as $N\to \infty$, where $d \nu_e(x) = \mu \bar{G}^s(x)dx$ and $d \eta_e(x)=\theta \bar{G}^r(x)dx$ for each $x \in \mathbb{R}_{+}$.
This, together with (\ref{eq:fluid-38}), (\ref{eq:fluid-39}), and $p_\opt = b_\opt \mu/\lambda$, gives that as $N\rightarrow\infty$,
\begin{align}
	&\bar{B}^N_\stat(0) =  \la 1,\bar{\nu}^N_\stat(0) \ra \Rightarrow \la 1, b_\opt \nu_e \ra  = b_{\opt}, \label{eq:thm1-aux1} \\
	&\bar{Q}^N_\stat(0) = \bar{X}^N_\stat(0) - \bar{B}^N_\stat(0) \Rightarrow b_\opt - b_\opt = 0. \label{eq:thm1-aux2}
\end{align}
The function $g_U$ is continuous. Hence, by~(\ref{eq:thm1-aux1}) and the continuous mapping theorem,
\begin{align}\label{eq:thm1-aux3}
	g_U(\bar{B}^N_\stat(0)) \Rightarrow g_U(b_{\opt}), \mbox{ as } N\rightarrow\infty.
\end{align}
Then, since $g_U$ is bounded, (\ref{eq:thm1-aux3}) and the bounded convergence theorem yield that
\begin{align}\label{eq:thm1-aux4}
	\lim_{N\to\infty}\mathbb{E}^N_{\xi}\left[g_U\left( \bar{B}_\stat^N(0) \right)\right] = g_U(b_{\opt}).
\end{align}
From (\ref{eq:Q-chi}) and (\ref{eq:thm1-aux2}),
\begin{align}\label{eq:thm1-aux5}
	\la 1_{[0,\chi^N_\stat(0)]}, \bar{\eta}^N_\stat(0) \ra = \bar{Q}^N_\stat(0) \Rightarrow 0, \mbox{ as } N\rightarrow\infty.
\end{align}
Note that for each $N \in \mathbb{N}$, 
\begin{align*}
	0 \leq a \la  1_{[0,\chi^N_\stat(0)]} h^r , \bar{\eta}^N_\stat(0) \ra \leq a \|h^r\|_{\infty} \bar{Q}^N_\stat(0),
\end{align*}
which, together with (\ref{eq:thm1-aux5}) and boundedness of $h^r$, implies 
\begin{align}\label{eq:thm1-aux6}
	a\la  1_{[0,\chi^N_\stat(0)]} h^r, \bar{\eta}^N_\stat(0) \ra \Rightarrow 0, \mbox{ as } N\rightarrow\infty.
\end{align}
By Lemma~\ref{lemma:stationary-exist}, $\lim_{N\to\infty}p_{\opt} \bar{\lambda}^N =p_{\opt} \lambda$, and \eqref{eq:thm1-aux0},
$$
\lim_{N\to\infty}{\mathbb E}_{\xi}^N\left[ \la 1, \bar{\eta}^N_\stat(0)\ra\right]=\lim_{N\to\infty}p_*\bar{\lambda}^N\theta^{-1}=p_*\lambda\theta^{-1}= \la 1, p_*\lambda\theta^{-1}\eta_e\ra.
$$
This together with \eqref{eq:thm1-aux0} implies that $\left\{\la 1, \bar{\eta}^N_\stat(0)\ra\right\}_{N\in\mathbb{N}}$ is uniformly integrable.
Note that $\la 1_{[0,\chi^N_\stat(0)]} h^r , \bar{\eta}^N_\stat(0) \ra \leq \|h^r\|_{\infty} \la 1, \bar{\eta}^N_\stat(0)\ra$ for each $N\in{\mathbb N}$ and $h^r$ is bounded. 
Thus,  $\left\{\la 1_{[0,\chi_\stat^N(0)]}h^r,\bar\eta_\stat^N(0)\ra\right\}_{N\in\mathbb{N}}$ is uniformly integrable.
This together with \eqref{eq:thm1-aux6} implies that
\begin{align}\label{eq:thm1-aux8}
	\lim_{N\to\infty} \mathbb{E}^N_{\xi}\left[a \la 1_{[0,\chi^N_\stat(0)]}h^r, \bar{\eta}^N_\stat \ra\right] = 0.
\end{align}
Finally, by Assumption~\ref{assumption:lambda-asymptotics}, it follows that 
\begin{align}\label{eq:thm1-aux9}
	\lim_{N\to\infty} a (1-p_{\opt}) \bar{\lambda}^N = a(1-p_{\opt})\lambda = a(\lambda-b_{\opt}\mu).
\end{align}
Combining (\ref{eq:thm1-aux4}), (\ref{eq:thm1-aux8}) and (\ref{eq:thm1-aux9}) establishes (\ref{eq:thm1-framework}), as desired. \qed\\


\textit{Proof of Theorem~\ref{theorem:lower-bound}. }
Fix a sequence $\{\hat{\pi}^N\}_{N\in\mathbb{N}}$ satisfying the conditions of Theorem~\ref{theorem:lower-bound}.
For each $N \in \mathbb{N}$, let $\xi^N\in{\mathcal S}(\hat{\pi}^N)$ be such that $\mathcal{C}^{N}(\hat{\pi}^{N})=\mathcal{C}_{\xi}^{N}(\hat{\pi}^{N})$ which exists by Lemma~\ref{lemma:stationary-exist} and the definition of $\mathcal{C}^{N}(\hat{\pi}^{N})$. For each $N \in \mathbb{N}$, let $Y^N_\stat$ be the state process for ($\hat{\pi}^N,\xi^N$).
It suffices to show that $\lim_{i \rightarrow\infty} \mathcal{C}_\xi^{N_i}(\hat{\pi}^{N_i}) \geq a(\lambda-b_\opt \mu) + g_U(b_\opt)$, for any convergent subsequence of cost functions $\left\{\mathcal{C}_\xi^{N_i}(\hat{\pi}^{N_i})\right\}_{i=1}^{\infty}$.  Fix such a subsequence $\{N_i\}_{i=1}^{\infty}$.
We consider the fluid scaled sequence $\{(\bar{X}^{N_i}_\stat, \bar{\nu}^{N_i}_\stat, \bar{\eta}^{N_i}_\stat)\}_{i=1}^{\infty}$. 
By Lemma~\ref{lemma:stationary-dist}, it suffices to show 
\begin{align}\label{eq:LB-aux-0}
	&\lim_{i \rightarrow\infty} \mathbb{E}^{N_i}_{\xi} \left[a  \left(1-\hat{p}^{N_i}\right)\bar{\lambda}^{N_i} + a  \la 1_{[0,\chi^{N_i}_\stat(0)]} h^r, \bar{\eta}^{N_i}_\stat(0) \ra + g_U\left(\bar{B}^{N_i}_\stat(0)\right)\right] \nonumber\\
	&\qquad \geq a(\lambda-b_\opt \mu) + g_U(b_\opt).
\end{align}

We begin by noting that the sequence $\{(\bar{X}^{N_i}_\stat(0), \bar{\nu}^{N_i}_\stat(0), \bar{\eta}^{N_i}_\stat(0))\}_{i=1}^{\infty}$ is tight. This follows by  Theorem~6.2 in~\cite{kang2012asymptotic} and its proof since in the present setting, the result in Lemma~6.1 in~\cite{kang2012asymptotic} holds, $\bar{K}^{N_i}_\stat(t) \leq \bar{E}^{N_i}_\stat(t) + \la 1, \bar{\eta}^{N_i}_\stat(0) \ra$ for all $i \in \mathbb{N}$ and $t \geq 0$, and $\bar{X}^{N_i}_\stat(0) \leq 1 + \la 1, \bar{\eta}^{N_i}_\stat(0) \ra$ for all $i\in\mathbb{N}$.   Since $\{(\bar{X}^{N_i}_{\stat}(0), \bar{\nu}^{N_i}_{\stat}(0), \bar{\eta}^{N_i}_{\stat}(0))\}_{i=1}^{\infty}$ is tight, there exists a further subsequence $\{N_{i_k}\}_{k=1}^{\infty}$ such that
\begin{align}\label{eq:LB-aux-1}
	\left(\bar{X}^{N_{i_k}}_\stat(0), \bar{\nu}^{N_{i_k}}_\stat(0), \bar{\eta}^{N_{i_k}}_\stat(0)\right) \Rightarrow \left(X^0_\stat, \nu^0_\stat,  \eta^0_\stat\right),
\end{align}
as $k\rightarrow\infty$. Without loss of generality, we can replace $\{N_{i_k}\}_{k=1}^{\infty}$ with $\{N_i\}_{i=1}^{\infty}$ by eliminating some members if necessary. In what follows, we verify that~(\ref{eq:LB-aux-0}) holds along this subsequence.  For this, we will first show that
\begin{align}
	&\lim_{i\to\infty}\mathbb{E}^{N_i}_{\xi} \left[ a\left(1-\hat{p}^{N_i}\right)\bar{\lambda}^{N_i} + a \la 1_{[0,\chi^{N_i}_\stat(0)]}h^r, \bar{\eta}^{N_i}_\stat(0) \ra + g_U\left(\bar{B}^{N_i}_\stat(0)\right) \right] \nonumber\\
	& =a\left(1-p\right)\lambda + a \mathbb{E}_{\xi} \left[\la  1_{[0,\chi_\stat^0]}h^r, \eta_\stat^0 \ra\right] + \mathbb{E}_{\xi} \left[g_U(B_\stat^0)\right],\label{eq:exp_conv}
\end{align}
where $\xi$ denotes the distribution of $\left(X^0_\stat, \nu^0_\stat,  \eta^0_\stat\right)$ and ${\mathbb E}_\xi$ is the expectation operator for $\xi$.  Then we will establish process level convergence to a stationary fluid model solution for $p\lambda$ in order to apply Lemma~\ref{lemma:thm1} to the right-hand side of \eqref{eq:exp_conv}.
	
We begin by showing that $\eta_\stat^0$ has no atoms and that $\left\{\la 1_{[0,\chi^{N_i}_\stat(0)]}h^r, \bar{\eta}^{N_i}_\stat(0) \ra\right\}_{i=1}^{\infty}$ is uniformly integrable.
By Lemma~\ref{lemma:stationary-exist}, Assumption~\ref{assumption:lambda-asymptotics} and $\lim_{i\to\infty}\hat p^{N_i}=p$,
\begin{equation}\label{eq:ui}
	\lim_{i\to\infty} \mathbb{E}^{N_i}_{\xi}\left[\la 1, \bar{\eta}^{N_i}_\stat(0)\ra\right]=p\lambda\theta^{-1},
\end{equation}
so $\sup_{i\in\mathbb{N}} \mathbb{E}^{N_i}_{\xi}\left[\la 1, \bar{\eta}^{N_i}_\stat(0)\ra\right]<\infty$.
This together with (\ref{eq:LB-aux-1}) implies that Assumption~\ref{assumption:asymptotics} holds with $\{E^{N_i}_{\hat{p}^{N_i}}\}_{i=1}^{\infty}$, $\{\hat{\pi}^{N_i}\}_{i=1}^{\infty}$ and $\{\xi^{N_i}\}_{i=1}^{\infty}$ replacing $\{E^N_{p^N}\}_{N\in\mathbb{N}}$, $\{\pi^N_{p^N}\}_{N\in\mathbb{N}}$ and $\{\varsigma^N\}_{N\in\mathbb{N}}$ respectively.
Thus, by Lemma~\ref{lemma:convergence-eta}, $\bar{\eta}^{N_i}_\stat \Rightarrow \eta_\stat$, as $i\to\infty$, where $\eta_\stat(0)\overset{d}{=}\eta_\stat^0$ and $\eta_\stat$ satisfies~(\ref{eq:fluid-47}) almost surely for $E_\stat(t)=p\lambda t$, $t \geq 0$.
Moreover, since $\bar{\eta}^{N_i}_\stat$ is a stationary process for each $i \in \mathbb{N}$, $\eta_\stat$ is a stationary process such that $\eta_\stat(t) \overset{d}{=} \eta_\stat^0$ for all $t \geq 0$.
Hence, by Lemma~\ref{lemma:etastat}, $\eta_\stat^0=p\lambda\theta^{-1}\eta_e$, so that $\eta_\stat^0$ has no atoms, $\la 1,\eta_\stat^0\ra=p\lambda\theta^{-1}$ and $x\mapsto\la 1_{[0,x]},\eta_\stat^0\ra$ is a continuous, strictly increasing function on ${\mathbb R}_+$. Then recalling \eqref{eq:ui}, $\lim_{i\to\infty}{\mathbb E}_{\xi}^{N_i}\left[ \langle 1,\bar\eta_\stat^{N_i}(0)\rangle\right]=\langle 1,\eta_\stat^0\rangle$.  This together with \eqref{eq:LB-aux-1} implies that $\left\{ \langle 1,\bar\eta_\stat^{N_i}(0)\rangle\right\}_{i=1}^{\infty}$ is uniformly integrable. Since $h^r$ is bounded and
$\la 1_{[0,\chi^{N_i}_\stat(0)]} h^r , \bar{\eta}^{N_i}_\stat(0) \ra \leq \|h^r\|_{\infty} \la 1, \bar{\eta}^{N_i}_\stat(0)\ra$ for each $i\in{\mathbb N}$,
uniform integrability of $\left\{\la 1_{[0,\chi^{N_i}_\stat(0)]}h^r, \bar{\eta}^{N_i}_\stat(0) \ra\right\}_{i=1}^{\infty}$ follows.

Next we show \eqref{eq:exp_conv}.
For this without loss of generality, we assume that the convergence in \eqref{eq:LB-aux-1} is almost sure which we abbreviate as a.s.   By \eqref{eq:LB-aux-1}, we have
\begin{align}\label{eq:convergence-B}
	&\lim_{i\to\infty}\bar{B}_\stat^{N_i}(0)=\lim_{i\to\infty}\langle 1,\bar\nu_\stat^{N_i}(0)\rangle=\langle 1,\nu_\stat^0\rangle=B_\stat^0,\quad\hbox{a.s., and} \\ &\lim_{i\to\infty}\bar Q_\stat^{N_i}(0)=\lim_{i\to\infty}\bar X_\stat^{N_i}(0)-\lim_{i\to\infty}\bar B_\stat^{N_i}(0)= X_\stat^0-B_\stat^0=Q_\stat^0,\quad\hbox{a.s.} \nonumber
\end{align}
This implies that
$$
\lim_{i\to\infty}\la 1_{[0,\chi^{N_i}_\stat(0)]}, \bar{\eta}^{N_i}_\stat(0)\ra=\lim_{i\to\infty}\bar Q_\stat^{N_i}(0)= Q_\stat^0=\la 1_{[0,\chi_\stat^0]}, \eta_\stat^0\ra,\quad\hbox{a.s.}
$$
Thus, since $x\mapsto\la 1_{[0,x]}, \eta_\stat^0 \ra$ is a continuous strictly increasing function on ${\mathbb R}_+$, $\lim_{i\to\infty}\chi^{N_i}_\stat(0)=\chi_\stat^0$ a.s.
This together with the above display and that $h^r$ is continuous and bounded implies
\begin{align}\label{eq:convergence-eta}
	\lim_{i\to\infty}\la 1_{[0,\chi^{N_i}_\stat(0)]}h^r, \bar{\eta}^{N_i}_\stat(0) \ra=
	\la 1_{[0,\chi_\stat^0]}h^r, \eta_\stat^0 \ra,\quad\hbox{a.s.}
\end{align}Now, as in the proof of Theorem~\ref{theorem:convergence-cost}, \eqref{eq:exp_conv} follows from  $\lim_{i\rightarrow\infty}\hat{p}^{N_i}\bar{\lambda}^{N_i}=p\lambda$, (\ref{eq:convergence-B}) and $g_U$ is bounded and continuous, and (\ref{eq:convergence-eta}) and the uniform integrability of $\left\{\la 1_{[0,\chi^{N_i}_\stat(0)]}h^r, \bar{\eta}^{N_i}_\stat(0) \ra\right\}_{i=1}^{\infty}$.
	
Finally, we argue process level convergence to a stationary fluid model solution for $p\lambda$.
Since Assumption~\ref{assumption:asymptotics} holds for $\{N_i\}_{i=1}^{\infty}$ (as noted above) and $\eta^0_\stat$ has no atoms (also noted above),
Lemma~\ref{lemma:weak-convergence} implies that $\left( \bar X_\stat^{N_i}, \bar \nu_\stat^{N_i}, \bar \eta_\stat^{N_i} \right) \Rightarrow \left( X_\stat, \nu_\stat, \eta_\stat \right)$, as $i\rightarrow\infty$,
where $\left( X_\stat, \nu_\stat, \eta_\stat \right)$ is almost surely a fluid model solution for $p\lambda$ such that $\left( X_\stat(0), \nu_\stat(0), \eta_\stat(0) \right)\overset{d}{=} \left(X^0_\stat, \nu^0_\stat,  \eta^0_\stat\right)$. Moreover, $\left(  X_\stat, \nu_\stat, \eta_\stat \right)$ is a stationary fluid model solution for $p\lambda$ by the stationarity of $(\bar{X}^{N_i}_\stat, \bar{\nu}^{N_i}_\stat, \bar{\eta}^{N_i}_\stat)$ for each $i \in \mathbb{N}$.  Then, from Lemma~\ref{lemma:thm1}, there exists $b\in[0,\min\{1,p\lambda/\mu\}]$ such that
$\mathbb{E}_\xi[B_\stat^0]=\mathbb{E}_\xi[B_\stat(0)]=b$. Since $g_U$ is convex, Jensen's inequality further implies that
\[
\mathbb{E}_{\xi} \left[ g_U(B_\stat^0) \right] \geq g_U \left( \mathbb{E}_{\xi} \left[B_\stat^0\right] \right) = g_U(b).
\]
This together with  \eqref{eq:exp_conv} and the second part of Lemma~\ref{lemma:thm1} gives 
\begin{align}\label{eq:thm-2-proof-last-step}
	&\lim_{i \rightarrow \infty} \mathbb{E}^{N_i}_{\xi}\left[ a(1-\hat{p}^{N_i})\bar{\lambda}^{N_i} + a\la 1_{[0,\chi^{N_i}_\stat(0)]} h^r, \bar{\eta}^{N_i}_\stat(0) \ra + g_U\left( \bar{B}^{N_i}_\stat(0)\right)\right] \nonumber \\
	\geq& a(1-p)\lambda +a\left(p\lambda- b \mu \right) + g_U \left( b \right) \nonumber\\
	=& a \left( \lambda - b \mu \right) + g_U \left( b \right) 
	\geq a \left( \lambda - b_{\opt} \mu \right) + g_U \left( b_{\opt} \right),
\end{align}
which completes the proof that \eqref{eq:LB-aux-0} holds, as desired.\qed

\vspace{-0.1in}
\appendix
\section{The Fluid Model for $\gamma$} \label{section:appendix}
\vspace{-0.1in}
We write the fluid model equations and write fluid model solutions for $\gamma>0$ in this appendix. We refer the reader to Section 3.1 in~\cite{puha2020fluid} for details. Given a Polish space $\mathbb{S}$, we use $\mathbf{C}(\mathbb{S})$ to denote the set of functions having domain $\mathbb{R}_{+}$ and range $\mathbb{S}$ that are continuous in time.

The fluid model for $\gamma$ has as an input a non-decreasing function $E(t)=\gamma t$, $t \geq 0$. We set $\mathbb{X}:=\mathbb{R}_{+} \times \mathbf{M}[0,H^s) \times \mathbf{M}[0,H^r)$, endowed with the product topology in a Polish space. To define the fluid model for $\gamma$, we consider $(X,\nu,\eta)\in\mathbf{C}(\mathbb{X})$ such that
\begin{align}
	\la 1_{\{x\}}, \eta(0) \ra = 0, \mbox{ for all } x\in[0,H^r),  \label{eq:fluid-no-atom}
\end{align}
and such that for each $t\geq 0$,
\begin{align}
	&\la 1,\nu(t) \ra \leq X(t) \leq \la 1,\nu(t) \ra + \la 1, \eta(t) \ra, \label{eq:fluid-35}\\
	& \la 1, \nu(t) \ra \leq 1, \label{eq:fluid-36}\\
	& \int_{0}^{t} \la h^s, \nu(u) \ra du < \infty \mbox{ and } \int_{0}^{t} \la h^r, \eta(u) \ra du <\infty. \label{eq:fluid-37}
\end{align}

Given $(X,\nu,\eta)\in\mathbf{C}(\mathbb{X})$ satisfying (\ref{eq:fluid-no-atom})-(\ref{eq:fluid-37})), we define auxiliary functions $B$, $Q$, $\chi$, $R$, $D$, and $K$ in $\mathbf{C}(\mathbb{R}_{+})$ and $I$ in $\mathbf{C}(\mathbb{R}_{+})$ as follows: for each $t\geq 0$,
\begin{align}
	&B(t) = \la 1, \nu(t) \ra, \label{eq:fluid-38}\\
	&Q(t) = X(t) - B(t), \label{eq:fluid-39}\\
	&\chi(t) = \inf\{x\geq 0: \la 1_{[0,x]}, \eta(t) \ra \geq Q(t)\}, \label{eq:def-chi} \\
	&R(t) = \int_{0}^{t} \left(\int_{0}^{\chi(u)} h^r(w)\eta(u)(dw)\right)du, \label{eq:fluid-40}\\
	&D(t) = \int_{0}^{t} \la h^s, \nu(u) \ra du, \label{eq:fluid-41}\\
	&K(t) = B(t) + D(t) - B(0), \label{eq:fluid-42}\\
	&I(t) = 1-B(t).\label{eq:fluid-43}
\end{align}
Then $B$, $Q$, $\chi$, $R$, $D$, $K$, and $I$ are fluid analogs of the busy server, the queue length, the waiting time of the HL fluid in queue, the reneging, the departure, the entry-into-service, and the idleness processes, respectively.

Further some additional properties and equations that should be satisfied by $(X,\nu,\eta)\in\mathbf{C}(\mathbb{X})$ are as follows: for any continuous and bounded function $f$ having domain $\mathbb{R}_{+}$, for each $t\geq 0$,
\begin{align}
	& K \mbox{ is non-decreasing}, \label{eq:fluid-44}\\
	& X(t) = X(0) + E(t) - R(t) - D(t), \label{eq:fluid-45}
\end{align}
\begin{align}
	& \la f, \nu(t) \ra = \la f(\cdot + t)\frac{\bar{G}^s(\cdot + t)}{\bar{G}^s(\cdot)}, \nu(0) \ra + \int_{0}^{t} f(t-u) \bar{G}^s(t-u) dK(u), \label{eq:fluid-46}\\
	& \la f, \eta(t) \ra = \la f(\cdot + t) \frac{\bar{G}^r(\cdot + t)}{\bar{G}^r(\cdot)}, \eta(0) \ra + \gamma\int_{0}^{t} f(t-u) \bar{G}^r(t-u) du. \label{eq:fluid-47}
\end{align}

\begin{definition}\label{definition:fluid-solution}
A fluid model solution for $\gamma>0$ is $(X,\nu,\eta)$ that satisfies (\ref{eq:fluid-no-atom})-(\ref{eq:fluid-37}), and~(\ref{eq:fluid-44})-(\ref{eq:fluid-47}).
\end{definition}

\begin{definition}\label{definition:fluid-solution-nonidling}
A non-idling fluid model solution for $\gamma>0$ is $(X,\nu,\eta)$ that satisfies Definition~\ref{definition:fluid-solution} and the following non-idling condition for each $t\geq 0$:
\begin{align} \label{eq:non-idling}
	I(t) = \left(1-X(t)\right)^+.
\end{align}	
\end{definition}









\vspace{-0.3in}
\bibliographystyle{ormsv080}
\bibliography{AsymptoticOptimalIdling_ZWP}

\section*{ONLINE APPENDIX}

\subsection{Proofs of Lemmas} \label{section:proofs}

Throughout this appendix, we fix $N\in{\mathbb N}$, $p\in(0,1]$, $\pi^N_p=({\mathbb S}^N,\{{\mathbb P}_y^N\}_{y\in{\mathbb S}^N})\in\Pi^N_p$ and a compatible initial distribution $\varsigma^N$, and we let $Y^N$ be the state process for $(\pi^N_p, \varsigma^N)$.
For each Borel subset $A$ of ${\mathbb S}^N$, define $L_0^N(A):=\mathbb{P}_{\varsigma}^N\left( Y^N(0)\in A\right)$ and for $t>0$, define
\[
L^N_t(A):=\frac{1}{t} \int_0^t \mathbb{P}_{\varsigma}^N\left( Y^N(s)\in A\right)ds.
\]
Then, for each $t>0$, $L_t^N$ is a probability measure on ${\mathbb S}^N$, and we use the notation $Y_t^N(0)=(\alpha_t^N(0),X_t^N(0),\nu_t^N(0),\eta_t^N(0))$ to denote a random vector with law $L_t^N$.
If $f:{\mathbb S}^N\to{\mathbb R}_+$ is measurable, then, for each $t\ge 0$, the expected value of $f(Y_t^N(0))$ is given by
$$
{\mathbb E}_{L_t}^N\left[ f(Y_t^N(0))\right]=\frac{1}{t}\int_0^t {\mathbb E}_{\varsigma}^N\left[ f(Y^N(s))\right]ds.
$$
Due to Lemma 4.4 in~\cite{kang2012asymptotic}, $\sup_{t\ge 0}{\mathbb E}_{\varsigma}^N\left[\la 1,\eta^N(t)\ra\right]<\infty$.  Thus, for all $t\ge 0$,
$$
{\mathbb E}_{L_t}^N\left[\la 1,\eta_t^N(0)\ra\right]<\infty,
$$
In addition, as shown in the proof of Proposition 4.1 in~\cite{atar2014fluid},
\begin{equation}\label{eq:LimExLt}
	\lim_{t\to\infty}{\mathbb E}_{L_t}^N\left[\la 1,\eta_t^N(0)\ra\right]=p\lambda^N\theta^{-1}.  
\end{equation}

\textit{Proof of Lemma~\ref{lemma:stationary-exist}. }
From Lemma 4.8 in~\cite{kang2012asymptotic} (which does not require the non-idling condition), the family of probability measures $\{ L^N_t \}_{t \geq 0}$ is tight. Since $\{ Y^N(t) \}_{t \geq 0}$ is a Feller Markov process such that ${\mathbb E}_{\varsigma}^N\left[\la 1,\bar{\eta}^N(0)\ra\right]<\infty$ (from Assumption~\ref{assumption:asymptotics}), the Krylov-Bogoliubov theorem (see Corollary 3.1.2 in~\cite{da1996ergodicity}) implies that any limit point $\xi^N$ of  $\{ L^N_t \}_{t \geq 0}$  is a stationary distribution for $\pi_p^N$ such that  if $Y_{\infty}^N(0)$ is a random variable with distribution $\xi^N$, then $\la 1,\eta_{\infty}^N(0)\ra$ has finite expected value under $\xi^N$.
Moreover, the marginal distribution of $\alpha_{\infty}^N(0)$ is such that the corresponding arrival process $E_{\infty}^N$ is a stationary renewal process with rate $p\lambda^N$.  Thus, $\la 1,\eta_{\infty}^N(0)\ra$  is equal in distribution to the stationary number of customers in a non-idling infinite server queue with arrival rate $p\lambda^N$ and service rate $\theta$.  Hence, for the admissible control policy $\pi^N_p$ with initial distribution $\xi^N$, we have ${\mathbb E}_{\xi}^N\left[\la 1,\eta_{\infty}^N(0)\ra \right]=p\lambda^N\theta^{-1}$ by Little's law. \qed \\


\textit{Proof of Lemma~\ref{lemma:stationary-dist}. }
Since $\pi^N_p$ is fixed, we suppress the process dependence on $\pi^N_p$ throughout the proof.
From the proof of Lemma~\ref{lemma:stationary-exist}, for each $N\in\mathbb{N}$, any limit point $\xi^N$ of  $\{ L^N_t \}_{t \geq 0}$  is a stationary distribution. Let $\{\tau(n)\}_{n\in\mathbb{N}}\subset\mathbb{R}_{+}$ be a strictly increasing subsequence along which $L_{\tau(n)}^N$  converges to $\xi^N$. On that subsequence,
\begin{align} \label{eq:utilization-limit}
	&\lim_{n \rightarrow \infty}\frac{1}{\tau(n)} \mathbb{E}_{\varsigma}^N \left[\int_0^{\tau(n)}  g_U\left( \bar{B}^N(t) \right) dt \right]\nonumber\\
	\overset{(1)}{=} & \lim_{n\rightarrow\infty} \frac{1}{\tau(n)} \int_{0}^{\tau(n)} \mathbb{E}^N_{\varsigma}\left[g_U\left(\bar{B}^N(t)\right)\right] dt \nonumber\\
	\overset{(2)}{=} & \lim_{n \rightarrow \infty} \mathbb{E}_{L_{\tau(n)}}^N\left[ g_U \left( \bar{B}^N_{\tau(n)}(0) \right) \right] \nonumber \\
	\overset{(3)}{=} & \mathbb{E}_\xi^N \left[g_U \left( \bar{B}^N_\stat(0) \right)\right].
\end{align}
where (1) follows by Fubini's theorem, (2) follows by definition of $L^N_{\tau(n)}$, and (3) follows because $g_U$ is continuous and bounded. Since $\{\tau(n)\}_{n\in\mathbb{N}}$ is arbitrary, (\ref{eq:utilization-limit}) implies (\ref{eq:long-run-avg-B}).

For each $t \geq 0$, let 
\[
M^N(t):=R^N(t)-\int_0^t \la 1_{[0,\chi^N(u-)]}h^r,\eta^N(u)\ra du.
\]
From Lemma 4 in~\cite{puha2020fluid}, $M^N(\cdot)$ is a martingale (with respect to the filtration $\{\mathcal{F}_t^N\}_{t\ge 0}$ defined in~\cite{puha2020fluid}). To see this, for $(x,u) \in [0,H^r) \times \mathbb{R}_{+}$, let $f^N(x,u) = 1_{[0,\chi^N(u-)]}(x)$ and note that $f^N$ is an almost surely bounded, measurable, real-valued function on $[0,H^r]\times \mathbb{R}_{+}$ so that Lemma~4 in~\cite{puha2020fluid} applies. Hence, for $n \in \mathbb{N}$,
\begin{align}\label{eq:R-Lt}
	\mathbb{E}_{\varsigma}^N\left[\frac{R^N(\tau(n))}{\tau(n)}\right]=&\mathbb{E}_{\varsigma}^N\left[\frac{1}{\tau(n)}\int_0^{\tau(n)}\la 1_{[0,\chi^N(u-)]}h^r,\eta^N(u)\ra du\right] \nonumber\\
	\overset{(2)}{=}&\mathbb{E}_{\varsigma}^N\left[\frac{1}{\tau(n)}\int_0^{\tau(n)}\la 1_{[0,\chi^N(u)]}h^r,\eta^N(u)\ra du\right] \nonumber\\
	\overset{(3)}{=}& \frac{1}{\tau(n)} \int_{0}^{\tau(n)} \mathbb{E}^N_{\varsigma}\left[\la 1_{[0,\chi^N(u)]}h^r,\eta^N(u)\ra\right]du \nonumber\\
	\overset{(4)}{=}&\mathbb{E}^N_{L_{\tau(n)}} \left[\la 1_{[0,\chi^N_{\tau(n)}(0)]} h^r, \eta^N_{\tau(n)}(0) \ra \right],
\end{align}
where (2) follows by noting that $\{ t\ge 0 : \chi^N(t-)\not=\chi^N(t) \}$ has Lebesgue measure zero, (3) follows by Fubini's theorem, and (4) follows by definition of $L^N_{\tau(n)}$.

By assumption, $\la 1,\eta_{\tau(n)}^N(0)\ra\Rightarrow\la 1,\eta_{\infty}^N(0)\ra$ as $n\to\infty$. From Lemma~\ref{lemma:stationary-exist}, ${\mathbb E}_\xi^N\left[\la 1,\eta_{\infty}^N(0)\ra\right]=p\lambda^N\theta^{-1}<\infty$. Hence, by \eqref{eq:LimExLt}, $\lim_{n\to\infty}{\mathbb E}_{L_{\tau(n)}}^N\left[\la 1,\eta_{\tau(n)}^N(0)\ra\right]={\mathbb E}_\xi^N\left[\la 1,\eta_{\infty}^N(0)\ra\right]$, implying that the sequence $\left\{\la 1,\eta_{\tau(n)}^N(0)\ra\right\}_{n\in{\mathbb N}}$ is uniformly integrable. Note that 
\begin{align}\label{eq:UI-inequality}
	\la 1_{[0,\chi^N_t(0)]}h^r,\eta^N_t(0)\ra \leq \Vert h^r\Vert_{\infty}\la 1,\eta^N_t(0) \ra, \mbox{ for all } t \geq 0.
\end{align}
Then, uniform integrability of $\left\{\la 1_{[0,\chi^N_{\tau(n)}(0)]}h^r,\eta^N_{\tau(n)}(0)\ra\right\}_{n\in{\mathbb N}}$ follows from~(\ref{eq:UI-inequality}), uniform integrability of $\left\{\la 1, \eta^N_{\tau(n)}(0)\ra\right\}_{N\in\mathbb{N}}$, and boundedness of $h^r$. Suppose we can show the following claim:
	\begin{claim}\label{claim:RenegingRateN} 
		$\la 1_{[0,\chi^N_{\tau(n)}(0)]} h^r, \eta^N_{\tau(n)}(0) \ra \Rightarrow \la 1_{[0,\chi^N_\stat(0)]} h^r, \eta^N_\stat(0) \ra$, as $n\to\infty$.
\end{claim}
Thus, (\ref{eq:R-Lt}), Claim~\ref{claim:RenegingRateN} and uniform integrability of $\left\{\la 1_{[0,\chi^N_{\tau(n)}(0)]}h^r,\eta^N_{\tau(n)}(0)\ra\right\}_{n\in{\mathbb N}}$ imply
\begin{align} \label{eq:reneging-limit}
	\lim_{n\to\infty}\mathbb{E}_{\varsigma}^N\left[\frac{R^N(\tau(n))}{\tau(n)}\right]
	=\mathbb{E}_{\xi}^N\left[\la
	1_{[0,\chi^N_\stat(0)]}h^r,\eta^N_\stat(0)\ra\right].
\end{align}
Then, since $\{\tau(n)\}_{n\in\mathbb{N}}$ is arbitrary, (\ref{eq:long-run-avg-R}) holds. \qed \\

To complete the proof, we verify Claim~\ref{claim:RenegingRateN} as follows.

\textit{Proof of Claim~\ref{claim:RenegingRateN}. } By assumption,
$\eta^N_{\tau(n)}(0) \Rightarrow \eta^N_\stat(0)$ as $n\to\infty$.  Without loss of generality, we may assume that this convergence is almost sure and we may fix an $w$ such that the stated convergence holds and evaluate all random elements at this $w$.
We have $\eta^N_{\tau(n)}(0)=\sum_{i=1}^{\la 1,\eta^N_{\tau(n)}(0)\ra}\delta_{w_i^n}$,  $\eta_\stat^N(0)=\sum_{i=1}^{\la 1,\eta^N_\stat(0)\ra}\delta_{w_i}$, and $\la 1,\eta^N_{\tau(n)}(0)\ra \Rightarrow \la 1,\eta_\stat^N(0)\ra$ as $n\to\infty$.  Since $\la 1,\eta^N_{\tau(n)}(0)\ra$ is non-negative integer valued, it follows that there exists a finite positive $\underline{n}$ such that for all $n\ge \underline{n}$, we have $\la 1,\eta^N_{\tau(n)}(0)\ra=\la 1,\eta_\stat^N(0)\ra$.  Hence, it follows that $w_i^n\to w_i$ for all $1\le i \le \la 1,\eta_\stat^N(0)\ra$, as $n\to\infty$.  Then, due to the continuity of $h^r$, as $n\to\infty$.
\begin{align*}
	\la 1_{[0,\chi^N_{\tau(n)}(0)]} h^r, \eta^N_{\tau(n)}(0) \ra&=\sum_{i=1}^{Q_{\tau(n)}^N(0)}h^r(w_i^n) \\
	&\to
	\sum_{i=1}^{Q_\stat^N(0)}h^r(w_i)=\la 1_{[0,\chi^N_\stat(0)]} h^r, \eta^N_\stat(0) \ra.
\end{align*}
\qed

\textit{Proof of Remark~\ref{remark:stationary-dist}. }
Let 
\begin{align*}
	\mathcal{C}^{N}_{\varsigma}(\pi^N_p,T) := \frac{1}{T} \mathbb{E}^{N}_{\varsigma} &\left[a \left( \bar{E}^N(T)-\bar{E}^N_p(T) + \bar{R}^N(T)\right) \right. \\
	& \left. + \int_{0}^{T} g_U\left(\bar{B}^N(t)\right)dt\right],
\end{align*}
	then $\mathcal{C}^N_{\varsigma}(\pi^N_p)=\limsup_{T\rightarrow\infty}\mathcal{C}^{N}_{\varsigma}(\pi^N_p,T)$.
	Let $\{\tau(n_i)\}_{i=1}^{\infty}$ be a subsequence of $\{\tau(n)\}_{n\in\mathbb{N}}$ such that $\{\mathcal{C}^{N}_{\varsigma}(\pi^N_p,\tau(n_i))\}_{i=1}^{\infty}$ converges and the limit is equal to the limit superior. Recalling that the sequence of probability measures $\{L_{\tau(n_i)}\}_{i=1}^{\infty}$ is tight, we denote by $\xi^N$ a limit point. Then, there exists a further subsequence $\{\tau(n_{i_k})\}_{k=1}^{\infty}$ such that $L^N_{\tau(n_{i_k})} \to \xi^N$ as $k\to\infty$. From (\ref{eq:utilization-limit}) and (\ref{eq:reneging-limit}), 
\begin{align*}
	&\mathcal{C}^N(\pi^N_p) \\
	=& \lim_{k\to\infty} \mathcal{C}^{N}_{\varsigma}(\pi^N_p,\tau(n_{i_k}))\\
	=&\lim_{k\to\infty}\mathbb{E}^N_{\varsigma} \left[a\frac{\bar{E}^N(\tau(n_{i_k}))-\bar{E}^N_p(\tau(n_{i_k}))}{\tau(n_{i_k})}\right] + \lim_{k\to\infty} \mathbb{E}^N_{\varsigma}\left[a\frac{\bar{R}^N(\tau(n_{i_k}))}{\tau(n_{i_k})}\right] \\
	&+\lim_{k\to\infty}\frac{1}{\tau(n_{i_k})}\mathbb{E}^N_{\varsigma} \left[\int_{0}^{\tau(n_{i_k})} g_U\left(\bar{B}^N(t)\right)dt\right]\\
	=&\mathbb{E}^N_{\xi} \left[a(1-p)\bar{\lambda}^N + a \la 1_{[0,\chi^N_\stat(0)]} h^r, \bar{\eta}^N_\stat(0) \ra + g_U\left(\bar{B}^N_\stat(0)\right)\right],
\end{align*}
which establishes the statement.\qed \\

\textit{Proof of Lemma~\ref{lemma:convergence-eta}.}
Due to Remark~\ref{remark:assumptions}, all but Assumptions~3(2) ($\lim_{N\to\infty}\mathbb{E}^N_{\xi}[\bar{X}^N(0)]=\mathbb{E}_{\xi}[X^0]<\infty$), 3(5) ($\lim_{N\to\infty}\mathbb{E}^N_{\xi}\left[\la 1, \bar{\eta}^N(0)\ra\right]=\mathbb{E}_{\xi}\left[\la 1, \eta^0\ra\right]<\infty$) and~5(2) ($\eta^0$ has no atoms) in~\cite{puha2020fluid} hold. A careful inspection of the proof of Theorem~1 in~\cite{puha2020fluid}, which relies on Theorem~6.2 in~\cite{kang2012asymptotic} (because the dynamics of $\eta^N$ are not altered by the non-idling condition), shows that $\sup_{N\in\mathbb{N}} \mathbb{E}^N_{\xi}\left[\la 1, \bar{\eta}^N(0) \ra\right] < \infty$ suffices in place of Assumptions 3(2) and 3(5) in~\cite{puha2020fluid}, and that Assumption~5(2) is not needed to establish that the limit $\eta$ exists and satisfies~(\ref{eq:fluid-47}). 
\qed \\

\textit{Proof of Lemma~\ref{lemma:etastat}.}
To see $\eta_\stat(t)=\gamma\theta^{-1}\eta_e$ for all $t \geq 0$ almost surely, note that $\eta_\stat$ satisfies \eqref{eq:fluid-47} almost surely. For each bounded continuous function $f$, the integrand of the first term of the right-hand side of \eqref{eq:fluid-47} tends to zero almost surely as $t\to\infty$ and so, by the dominated convergence theorem, the integral also tends to zero almost surely as $t\to\infty$.  In addition, $E_\stat(u)=\gamma u$ for all $u\ge 0$ and so the second term of the right-hand side of  \eqref{eq:fluid-47} converges to $\gamma\int_0^{\infty} f(u)(1-G^r(u))du=\gamma\theta^{-1}\la f,\eta_e\ra$ almost surely, as $t\to\infty$.  Thus, $\lim_{t\to\infty}\langle f,\eta_\stat(t)\rangle=\gamma\theta^{-1}\la f,\eta_e\ra$ for each bounded continuous function $f$. Since $\eta_\stat$ is a stationary process, the result follows. 
\qed\\

\textit{Proof of Lemma~\ref{lemma:thm1}.}
For $t\ge 0$, we have
\begin{align}\label{eq:claim-K}
\mathbb{E}_{\xi}[K_\stat(t)] \overset{(1)}{=} \mathbb{E}_{\xi}[D_\stat(t)] &\overset{(2)}{=} \mathbb{E}_{\xi}\left[\int_{0}^{t} \la h^s, \nu_\stat(u) \ra du \right] \nonumber\\
&\overset{(3)}{=} \mathbb{E}_{\xi}[\la h^s, \nu_\stat(0)\ra]t,
\end{align}
where (1) follows from (\ref{eq:fluid-42}) and $\mathbb{E}_{\xi}[B_\stat(t)]=\mathbb{E}_{\xi}[B_\stat(0)]$ because $B_\stat$ is a stationary process, (2) follows from (\ref{eq:fluid-41}), and (3) follows from Fubini's theorem and the stationarity of $\nu_\stat$.

From (\ref{eq:fluid-46}), for each $t\ge 0$,
\begin{align*}
\la 1, \nu_\stat(t) \ra =& \int_{0}^{\infty} \frac{\bar{G}^s(x+t)}{\bar{G}^s(x)} \nu_\stat(0)(dx) + \int_{0}^{t}\bar{G}^s(t-u) dK_\stat(u).
\end{align*}
Taking expectation on both sides of the above display and noting that $\la 1, \nu_\stat(t)\ra=B_\stat(t)$ from (\ref{eq:fluid-38}) yields that for each $t\ge 0$
\begin{align*}
\mathbb{E}_{\xi}[B_\stat(t)]
=& \mathbb{E}_{\xi}\left[\int_{0}^{\infty} \frac{\bar{G}^s(x+t)}{\bar{G}^s(x)} \nu_\stat(0)(dx)\right] \\
+& \mathbb{E}_{\xi}\left[\int_{0}^{t}(\bar{G}^s(t-u)) dK_\stat(u)\right].
\end{align*}
As $t\to\infty$, the first expectation converges to zero by two applications of the dominated convergence theorem and the fact that the interior integrand converges to zero almost surely as $t\to\infty$. Hence, since the left-hand side of the above is constant in $t$, combining this with (\ref{eq:claim-K}) gives that, for all $t\geq 0$,
\begin{align}
\mathbb{E}_{\xi}[B_\stat(t)]=&\mathbb{E}_{\xi}[\la h^s, \nu_\stat(0)\ra] \int_{0}^{\infty}(1-G^s(u))du\nonumber  \\
=&\mathbb{E}_{\xi}[\la h^s, \nu_\stat(0)\ra]\mu^{-1}=\mathbb{E}_{\xi}[D_\stat(t)](t\mu)^{-1}.\label{eq:bstat}
\end{align}
Furthermore, since $X_\stat(t) = X_\stat(0)+E_\stat(t)-R_\stat(t)-D_\stat(t)$, $t\ge 0$, (from (\ref{eq:fluid-45}) and  $\mathbb{E}_{\xi}[X_\stat(t)]=\mathbb{E}_{\xi}[X_\stat(0)]$, $t\ge 0$, due to the stationarity of $X_\stat$, it follows that $\mathbb{E}_{\xi}[D_\stat(t)] = \mathbb{E}_{\xi}[E_\stat(t)] - \mathbb{E}_{\xi}[R_{\stat}(t)] \leq \mathbb{E}_{\xi}[E_\stat(t)]=\gamma t$ for all $t\geq 0$. Then, \eqref{eq:bstat} implies that $\mathbb{E}_{\xi}[B_\stat(t)]\leq \gamma/\mu$, $t\ge 0$. Also, $B_\stat(t)\le 1$, $t\ge 0$.  Hence, there exists $b\in[0,\min\{1,\gamma/\mu\}]$ such that $\mathbb{E}_{\xi}[B_\stat(t)]=b$ for all $t\geq 0$. This together with~\eqref{eq:bstat} gives that, for all $t\ge 0$,
\begin{align*}
\mathbb{E}_{\xi}[D_\stat(t)] = b\mu t,
\end{align*}
which implies that for all $t\ge 0$,
\begin{align*}
\mathbb{E}_{\xi}[R_\stat(t)] = \mathbb{E}_{\xi}[E_p(t)] - \mathbb{E}_{\xi}[D_\stat(t)] = (\gamma-b\mu) t.
\end{align*}
From (\ref{eq:fluid-40}), Fubini's theorem and stationarity, for all $t\ge 0$,
\begin{align*}
\mathbb{E}_{\xi}[R_\stat(t)] =& \mathbb{E}_{\xi}\left[\int_{0}^{t} \la  1_{[0,\chi_\stat(u)]}h^r, \eta_\stat(u) \ra du\right] \\
=& \mathbb{E}_{\xi}\left[\la  1_{[0,\chi_\stat(0)]}h^r, \eta_\stat(0) \ra\right]t.
\end{align*}
The above two displays imply that $\mathbb{E}_{\xi}\left[\la  1_{[0,\chi_\stat(t)]}h^r, \eta_\stat(t) \ra\right]=\gamma-b\mu$ for all $t \geq 0$. 
\qed\\

\subsection{Extension: Holding Cost} \label{section:extension}
In this appendix, we include holding costs in the objective function to penalize congestion, and we show similar results as in the case with abandonment cost only, using the enlarged admissible policy class with admission control.

Let $c$ be the holding cost incurred per customer per unit time. Then, given $\pi^N \in \Pi^N$ (Definition~\ref{definition:admissible-policy}) and a compatible initial distribution $\varsigma^N$, the long-run average cost of $(\pi^N, \varsigma^N)$ is modified as
\begin{align*}
	\mathcal{H}_{\varsigma}^N(\pi^{N}) := \limsup_{T\rightarrow\infty} \frac{1}{T} \mathbb{E}^{N}_{\varsigma} &\left[c\int_{0}^{T}\bar{Q}^{N}(\pi^{N},t)dt +
	a \bar{R}^{N}(\pi^{N},T) \right. \\
	& \left. + \int_{0}^{T} g_U\left(\bar{B}^{N}( \pi^{N},t)\right)dt\right],
\end{align*}
and the worst case cost under $\pi^N$ is
\begin{align*}
	\mathcal{H}^N(\pi^N) := \sup_{\xi^N\in\mathcal{S}(\pi^N)} \mathcal{H}^N_{\xi}(\pi^N).
\end{align*}
Then, the objective is to determine $\pi_{\mbox{\tiny opt}}^{N}$ and $\mathcal{H}^N(\pi_{\mbox{\tiny opt}}^{N})$ such that
\begin{align} \label{eq:scheduling-generalize}
	\mathcal{H}^N(\pi_{\mbox{\tiny opt}}^{N}) := \inf_{\pi^N\in\Pi^N} \mathcal{H}^N(\pi^N).
\end{align}

We begin by noting that the fluid model equations are not changed, and hence the fluid control problem can be obtained based on the unchanged fluid model invariant states. Also, the weak convergence result (Lemma~\ref{lemma:weak-convergence}) continues to hold because the proof of Lemma~\ref{lemma:weak-convergence} does not rely on the objective function.

\begin{assumption} \label{assumption:hDFR}
	The function $h^r$ is non-increasing.
\end{assumption}
Assumption~\ref{assumption:hDFR} is crucial to prove the main asymptotic optimality result, Theorem~\ref{theorem:lower-bound}, when holding costs are considered. To explain this point, for $p\in(0,1]$ and $b\in[0,\min\{1,p\lambda/\mu\}]$, define
\begin{align} \label{eq:def-q}
	q(b,p) := p\lambda \int_{0}^{(G^r)^{-1}(1-b\mu/p\lambda)}(1-G^r(x)) dx.
\end{align}
From Equation~(54) in~\cite{puha2020fluid}, $q(b,p)$ represents the invariant fluid queue length for $p\in(0,1]$ and $b\in[0,\min\{1,p\lambda/\mu\}]$ when the arrival rate is thinned to $p\lambda$. As a consequence of Assumption~\ref{assumption:hDFR}, $q(\cdot,p)$ is a convex function on $[0,\min\{1,p\lambda/\mu\}]$ for each $p\in(0,1]$ (see Remark~10 in~\cite{puha2019scheduling}). This convexity plays an integral role in our analysis. The modified fluid control problem is
\begin{align} \label{eq:scheduling-fluid-generalize}
	\min_{b\in[0,\min\{1,\lambda/\mu\}]} \  c  q(b,1) + a(\lambda-b\mu) + g_U(b).
\end{align}
Different from~\eqref{eq:scheduling-fluid} when the holding costs were not included, the optimization problem~\eqref{eq:scheduling-fluid-generalize} becomes sensitive to the patience distribution, because the fluid queue length depends on the patience distribution; see the right-hand side of (\ref{eq:def-q}).
We denote the solution to~(\ref{eq:scheduling-fluid-generalize}) by $b_{\opt,\hc}$ (which is unique and guaranteed to exist because $q(b,1)$ is a convex function). As in Section~\ref{section:fluid}, if $b_{\opt,\hc} < \min\{1,\lambda/\mu\}$, then we expect an idling control policy to be optimal for~\eqref{eq:scheduling-generalize}. In particular,  servers can be allowed to take a rest for $(1-b_{\opt,\hc})(b_{\opt,\hc}\mu)^{-1}$ time units after each service completion.

As in Sections~\ref{section:policy} and~\ref{section:ao}, in order to show the asymptotic optimality property, we work with the enlarged admissible policy class $\hat{\Pi}^N$, which incorporates the potential of admission control. The unit abandonment and holding cost for the admitted arrivals remain $a$ and $c$. For every rejected arrival, in addition to a cost of $a$ (as in~\eqref{eq:scheduling-fluid-2}), we need to further account for the holding cost that would have been incurred if under a control policy in $\Pi^N$ without admission control that may idle. 
Suppose that $\tilde{\mathcal{H}}(b,p)$ is the overall holding costs for the rejected arrivals for $p\in(0,1]$ and  $b\in[0,\min\{1,p\lambda/\mu\}]$, and we call it the fluid-scaled holding cost compensator.

\begin{definition}[The Fluid-Scaled Holding Cost Compensator]\label{definition:hc-compensator}
	Given $\lambda$ is fixed, $p\in(0,1]$ and  $b\in[0,\min\{1,p\lambda/\mu\}]$, $\tilde{\mathcal{H}}(b,p)$ is given by
	\begin{align*}
		\tilde{\mathcal{H}}(b,p)= c\left(q(b,1)-q(b,p)\right).
	\end{align*}	
\end{definition}

Then, the modified fluid control problem under $\hat{\Pi}^N$ is
\begin{align} \label{eq:scheduling-fluid-generalize-2}
	\min_{p\in(0,1], b\in[0,\min\{1,p\lambda/\mu\}]} \ &c q(b,p) + \tilde{\mathcal{H}}(b,p)+ a(p\lambda-b\mu) \nonumber \\
	&  + a(1-p)\lambda + g_U(b)\nonumber \\
	=\min_{b\in[0,\min\{1,p\lambda/\mu\}]}\ &c q(b,1) + a(\lambda-b\mu) + g_U(b).
\end{align}
Under Definition~\ref{definition:hc-compensator}, it is clear that the solution to~\eqref{eq:scheduling-fluid-generalize-2} is identical to the solution to~\eqref{eq:scheduling-fluid-generalize}. Then, from~\eqref{eq:scheduling-fluid-generalize-2}, given $p \in (0,1]$, $\pi^{N}_p\in\Pi^N_p$ (Definition~\ref{definition:admissible-policy-2}) and a compatible initial distribution $\varsigma^N$, the modified objective function of $(\pi^N_p, \varsigma^N)$ is given by
\begin{align}\label{eq:scheduling-generalize-2}
	&\mathcal{H}_{\varsigma}^N(\pi^{N}_p) = \limsup_{T\rightarrow\infty} \frac{1}{T} \mathbb{E}^{N}_{\varsigma} \left[c\int_{0}^{T} \bar{Q}^N(\pi^{N}_p,t) dt \right.  \nonumber\\
	& +  \int_{0}^{T}\tilde{\mathcal{H}}\left(\bar{B}^N(\pi^{N}_p,t),p\right)dt + a \left(\bar{E}^N(T)-\bar{E}^N_p(T)+\bar{R}^{N}(\pi^{N}_p,T)\right)  \nonumber\\ 
	& \left. + \int_{0}^{T} g_U\left(\bar{B}^N(\pi^{N}_p,t)\right)dt\right].
\end{align}

Noting that the fluid solution to~\eqref{eq:scheduling-fluid-generalize}, $b_{\opt,\hc}$, is independent of $p$, we can consider the same proposed policy as defined in Definition~\ref{definition:policy} with $b_{\opt}$ replaced by $b_{\opt,\hc}$; denote it by $\pi^N_{\opt,\hc}$. In the remainder of the appendix, we outline the proof of asymptotic optimality for $\{\pi^N_{\opt,\hc}\}_{N\in\mathbb{N}}$.

Lemma~\ref{lemma:stationary-exist} continues to hold, because its proof does not rely on the objective function. For the objective~(\ref{eq:scheduling-fluid-generalize}) with linear holding cost penalties, Lemma~\ref{lemma:stationary-dist} can be modified such that for any admissible HL control policy $\pi^N_p \in \Pi^N_p$, the following holds for a stationary distribution $\xi^N\in\mathcal{S}(\pi^N_p)$:
\begin{align*}
	\mathcal{H}^N_{\varsigma}(\pi^N_p) = \mathbb{E}^N_{\xi}& \Bigg[c \bar{Q}^N_\stat(0) + \tilde{\mathcal{H}}\left(\bar{B}^N_\stat(0),p\right)	\\
	+& a(1-p)\bar{\lambda}^N+ a\la 1_{[0,\chi^N_\stat(0)]} h^r, \bar{\eta}^N_\stat(0) \ra + g_U\left(\bar{B}^N_\stat(0)\right)\Bigg].
\end{align*}
This is because for each $N\in\mathbb{N}$, if $\{\tau(n)\}_{n\in\mathbb{N}}\subset\mathbb{R}_{+}$ is a strictly increasing subsequence along which $L^N_{\tau(n)}$ converges to $\xi^N$, then $Q^N_{\tau(n)}(\pi^N_p,0) \Rightarrow Q^N_\infty(0)$ as $n\to\infty$, and
\begin{align*}
	\lim_{n \rightarrow \infty}\frac{1}{\tau(n)} \mathbb{E}_{\varsigma}^N \left[\int_0^{\tau(n)} Q^N(\pi^N_p,t) dt\right] =& \lim_{n\to\infty} \mathbb{E}^N_{L_{\tau(n)}} [Q^N_{\tau(n)}(0)] \\
	=& \mathbb{E}^N_{\xi}\left[Q^N_{\stat}(0)\right],
\end{align*}
because $Q^N(t) \leq \la 1,\eta^N(t)\ra$ for all $t\geq 0$, and $\left\{\la 1, \eta^N_{\tau(n)}(0) \ra\right\}_{n\in\mathbb{N}}$ is uniformly integrable as in the proof of Lemma~\ref{lemma:stationary-dist} in the online appendix~\ref{section:proofs} (see the paragraph immediately above Claim~\ref{claim:RenegingRateN}).

Thus, following a similar proof to that of Theorem~\ref{theorem:convergence-cost}, the sequence $\{\pi^N_{\opt,\hc}\}_{N\in\mathbb{N}}$ satisfies
\begin{align*}
	&\lim_{N\rightarrow\infty} \mathcal{H}^N(\pi^N_{\opt,\hc}) \\
	=& c q\left(b_{\opt,\hc},p\right) + \tilde{\mathcal{H}}(b_{\opt,\hc},p)+ a\left( \lambda - b_{\opt,\hc} \mu \right) +  g_U(b_{\opt,\hc})\\
	=& c q(b_{\opt,\hc},1)+ a\left( \lambda - b_{\opt,\hc} \mu \right) +  g_U(b_{\opt,\hc}).
\end{align*}
We wish to show that for $\{\hat{\pi}^N\}_{N\in\mathbb{N}}$ that satisfy the conditions in the statement of Theorem~\ref{theorem:lower-bound}, 
\begin{align*}
	\liminf_{N\rightarrow\infty} \mathcal{H}^N(\hat{\pi}^N) \geq c q(b_{\opt,\hc},1) + a\left( \lambda - b_{\opt,\hc} \mu \right) + g_U(b_{\opt,\hc}). 
\end{align*}

Let $\{( \bar X_\stat^{N_i}, \bar \nu_\stat^{N_i}, \bar \eta_\stat^{N_i} )\}_{N\in\mathbb{N}}$ and $\left( X_\stat, \nu_\stat, \eta_\stat \right)$ be as in the proof of Theorem~\ref{theorem:lower-bound}; that is,
\[
( \bar X_\stat^{N_i}, \bar \nu_\stat^{N_i}, \bar \eta_\stat^{N_i} ) \Rightarrow \left( X_\stat, \nu_\stat, \eta_\stat \right), \mbox{ as } i \rightarrow \infty,
\]
with  $\left( X_\stat, \nu_\stat, \eta_\stat \right)$ being almost surely a stationary fluid model solution for $p\lambda$. Note that, by definition,
\begin{align*}
	&c  Q_\stat(0) + \tilde{\mathcal{H}}(p,B_\stat(0)) \\
	=& c q(B_\stat(0),p) + c\left(q(B_\stat(0),1)-q(B_\stat(0),p)\right) \\
	=& c q(B_\stat(0),1).
\end{align*}
Since $q(\cdot,1)$ is convex on $[0,\min\{1,\lambda/\mu\}]$, Jensen's inequality together with Lemma~\ref{lemma:thm1} implies that
\begin{align}\label{eq:hc-convex}
	&\mathbb{E}_{\xi}[c Q_\stat(0) + \tilde{H}(p,B_\stat(0))]
	= \mathbb{E}_{\xi}[c q(B_\stat(0),1)] \nonumber\\
	\geq& c q(\mathbb{E}_{\xi}[B_\stat(0)],1) = c q(\mathbb{E}_{\xi}[B_\stat^0],1) = c q(b,1).
\end{align}

Then, as in the proof of Theorem~\ref{theorem:lower-bound},
\begin{align*}
	&\lim_{i \rightarrow \infty} \mathbb{E}^{N_i}_{\xi}\left[c  \bar{Q}^{N_i}_\stat(0) + \tilde{\mathcal{H}}\left(\bar{B}^{N_i}_\stat(0),p\right) \right.\\
	& \quad\quad\quad + a(1-\hat{p}^{N_i})\bar{\lambda}^{N_i} + \left. a\la 1_{[0,\chi^{N_i}_\stat(0)]} h^r, \bar{\eta}^{N_i}_\stat(0) \ra + g_U\left(\bar{B}^{N_i}_\stat(0) \right)\right]\\
	=& a(1-p)\lambda + \mathbb{E}_{\xi}\left[c Q_\stat(0) + \tilde{\mathcal{H}}(B_\stat(0),p) \right.\\
	&\quad\quad\quad \left. + a\la 1_{[0,\chi_\stat(0)]} h^r, \eta_\stat(0)\ra + g_U(B_\stat(0))\right]\\
	\geq& c q(b,1) + a(\lambda-b\mu) + g_U(b)\\
	\geq& c q(b_{\opt,\hc},1) + a\left( \lambda - b_{\opt,\hc} \mu \right) + g_U(b_{\opt,\hc}),
\end{align*}
where the first inequality follows from the first inequality in~(\ref{eq:thm-2-proof-last-step}) in the proof of Theorem~\ref{theorem:lower-bound}, and (\ref{eq:hc-convex}). 

Hence, the proposed policy $\pi^N_{\opt,\hc}$ with $p_{\opt,\hc}=b_{\opt,\hc}\mu/\lambda$ is asymptotically optimal, under the generalized objective function involving holding costs.

\begin{remark}
All the results in this section continue to stand if we consider a non-decreasing, continuous, convex holding cost that maintains uniform integrability of the sequence of holding costs rather than a linear holding cost. In that case, the corrected fluid-scaled holding cost compensator defined in Definition~\ref{definition:hc-compensator} is the difference of the holding cost function evaluated at $q(b,1)$ and at $q(b,p)$.
\end{remark}

\end{document}